\documentclass[11pt]{amsart}
\usepackage[a4paper,left=3.5cm,right=3.5cm,top=3.5cm,bottom=3.5cm]{geometry}

\makeatletter
\def\subsection{\@startsection{subsection}{2}
  \z@{.5\linespacing\@plus.7\linespacing}{.5\linespacing}
  {\normalfont\bfseries}}
\makeatother

\makeatletter
\def\subsubsection{\@startsection{subsubsection}{3}
  \z@{.4\linespacing\@plus.6\linespacing}{.4\linespacing}
  {\normalfont\mdseries}}
\makeatother

\usepackage[colorlinks=true, linkcolor=blue, citecolor=red, backref=page]{hyperref}
\renewcommand*{\backrefalt}[4]{
  \ifcase #1
    \relax
  \or
    {\small $\uparrow$ p.~#2}
  \else
    {\small $\uparrow$ pp.~#2}
  \fi
}

\usepackage{upref}
\usepackage{amssymb}
\usepackage{amsxtra}
\usepackage{mathtools}
\usepackage{tikz}
\usetikzlibrary{calc, patterns, intersections, decorations.pathreplacing}
\usepackage{multirow}
\usepackage{mathrsfs}

\newtheorem{Thm}{Theorem}[section]
\newtheorem{Lem}[Thm]{Lemma}
\newtheorem{Prop}[Thm]{Proposition}
\newtheorem{Cor}[Thm]{Corollary}

\newtheorem{Conj}[Thm]{Conjecture}

\theoremstyle{definition}
\newtheorem{Def}[Thm]{Definition}
\newtheorem{Ex}[Thm]{Example}

\theoremstyle{remark}
\newtheorem{Rem}[Thm]{Remark}

\numberwithin{equation}{section}
\allowdisplaybreaks[4]

\newcommand{\dc}{\overset{\raisebox{-0.4ex}{\scalebox{0.5}{$\circ$}}}{\cap}}

\title[Boundedness of Lusztig's $\boldsymbol{a}$-function]{Boundedness of Lusztig's $\boldsymbol{a}$-function for Coxeter groups of finite rank}

\author{Xiaoyu Chen}
\address{(Xiaoyu Chen) \newline \indent Department of Mathematics, Shanghai Normal University, 100 Guilin Road, Shanghai 200234, P.\ R.\ China}
\email{chenxiaoyu@shnu.edu.cn}

\author{Hongsheng Hu}
\address{(Hongsheng Hu) \newline \indent School of Mathematics, Hunan University, Changsha 410082, P.\ R.\ China}
\email{huhongsheng@hnu.edu.cn}
\date{September 10, 2025}



\begin{document}

\begin{abstract}
We prove that there exists a bound $N'_L(W)$ for a positively weighted Coxeter group $(W, S, L)$ of finite rank.
In particular, Lusztig's $\boldsymbol{a}$-function of $(W, S, L)$ is bounded.
\end{abstract}
\maketitle

\setcounter{tocdepth}{2}
\tableofcontents

\section{Introduction}

\subsection{Overview}

In his seminal paper \cite{Lus85}, Lusztig introduced the $\boldsymbol{a}$-function on a Coxeter group $(W,S)$, which has proven to be an important tool in the study of Kazhdan--Lusztig cells and many representation-theoretic topics.
For affine Weyl groups, Lusztig proved in \cite{Lus85} that the $\boldsymbol{a}$-function is bounded by the length of the longest element of the corresponding Weyl group.
For general Coxeter groups of finite rank, Xi \cite[Question~1.13(iv)]{Xi94} conjectured that the $\boldsymbol{a}$-function is bounded by the maximal length of elements of all the finite standard parabolic subgroups of $W$.

More generally, in Lusztig's monograph \cite{Lus03}, the boundedness conjecture was stated for a positively weighted Coxeter group $(W, S, L)$ of finite rank.
To be precise, the generator set $S$ is assumed to be finite and $L: W \to \mathbb{N}$ is a function (called a weight function) such that $L(s) > 0$ for any $s\in S$ and $L(ww') = L(w) +L(w')$ for any $w, w' \in W$ such that $\ell(ww') = \ell(w) + \ell(w')$, where $\ell : W \to \mathbb{N}$ is the length function.
Let $v$ be an indeterminate and $\mathcal{A} :=\mathbb{Z}[v^{\pm 1}]$ be the ring of Laurent polynomials in $v$ with integer coefficients.
Following \cite{Lus03}, the associated Hecke algebra $\mathcal{H}$ is a free $\mathcal{A}$-module with basis $\{T_w \mid w\in W\}$ and the multiplication is given by
\begin{equation*}
\begin{cases}
T_w T_{w'}=T_{ww'}, & \text{ if } \ell(ww')=\ell(w)+\ell(w'),\\
T_s^2= (v^{L(s)} - v^{-L(s)}) T_s + T_e, & \text{ if } s\in S,
\end{cases}
\end{equation*}
where $e$ is the neutral element of $W$. For any $x,y\in W$, we write
\begin{equation}\label{eq-strconst}
T_x T_y=\sum_{z\in W} f_{x,y,z} T_z, \quad f_{x,y,z}\in\mathcal{A}.
\end{equation}
We say that $N \in \mathbb{N}$ is a bound for $(W,S,L)$ if $v^{-N} f_{x,y,z} \in \mathbb{Z}[v^{-1}]$ for all $x, y, z \in W$, and that $(W,S,L)$ is bounded if there exists a bound for $(W,S,L)$.
We denote
\[N_L(W) := \max_I \{L(w_I)\}\]
where $I$ runs over subsets of $S$ such that the standard parabolic subgroup $W_I$ generated by $I$ is finite, and $w_I$ is the longest element in $W_I$.
The following conjecture is a generalization of the above mentioned \cite[Question~1.13(iv)]{Xi94} (and equivalent in the case $L = \ell$).

\begin{Conj} [{\cite[Conjecture~13.4]{Lus03}}] \label{conj-bound}
  The number $N_L(W)$ is a bound for any positively weighted Coxeter group $(W, S, L)$ of finite rank.
\end{Conj}

\begin{Rem}
    In his note \cite[Problem~(b)]{Lus20}, Lusztig also conjectured the existence of a bound for $(W,S,L)$ but did not specify the explicit bound $N_L(W)$.
\end{Rem}

Clearly, $(W, S, L)$ is bounded if $W$ is finite.
Moreover, $N_L(W)$ is a bound in this case; see \cite[Lemma~13.3]{Lus03}.
For infinite $W$, Conjecture~\ref{conj-bound} has been verified in the following cases:
\begin{enumerate}
    \item in \cite{Lus85}, Lusztig proved the boundedness for affine Weyl groups in the case $L = \ell$, and he pointed out (in the remark just below \cite[Conjecture~13.4]{Lus03}) that the same proof remains valid for weighted affine case;
    \item in \cite{Belolipetsky}, Belolipetsky proved the conjecture in the case $L = \ell$ for right-angled Coxeter groups, that is, the order of $st$ is either 2 or $\infty$ for any $s,t\in S$ with $s \ne t$;
    \item in \cite{Xi12}, Xi proved the conjecture in the case $L = \ell$ for $(W,S)$ with complete Coxeter graph;
    \item in \cite{Z}, Zhou proved the conjecture in the case $L = \ell$ for $(W,S)$ of rank 3;
    \item in \cite{SY15}, Shi and Yang proved the conjecture for weighted universal Coxeter groups;
    \item in \cite{SY}, Shi and Yang proved the conjecture for weighted Coxeter groups with complete Coxeter graph;
    \item in \cite{LS}, Li and Shi proved the conjecture for weighted Coxeter groups in which the order of $st$ is not $3$ for any $s,t\in S$;
    \item in \cite{Gao22}, Gao proved the conjecture for weighted Coxeter groups of rank $3$.
\end{enumerate}

The boundedness property for $(W,S,L)$ is a foundation of many research related to the structure of Kazhdan--Lusztig cells and cell representations.
Here are some instances.
\begin{enumerate}
  \item In \cite[Section~14]{Lus03}, for each $(W, S, L)$ Lusztig proposed his conjectures P1--P15 built on the assumption that $(W, S, L)$ is bounded (without assuming that $N_L(W)$ is a bound).
      In the case $L = \ell$, assuming $(W, S, \ell)$ is bounded, these conjectures are verified; see \cite[Section~15]{Lus03} and \cite{Lus87}.
      The proof utilizes the Kazhdan--Lusztig positivity conjecture which was proved in the landmark paper \cite{EW}.
      More generally, conjectures P1--P15 are verified in the so-called quasisplit case with the boundedness assumption; see \cite[Section~16]{Lus03}.
  \item In some other cases where the boundedness in Conjecture~\ref{conj-bound} is established, Lusztig's conjectures P1--P15 are proved; for example, \cite{GX21, Geck11, GP19-G2, SY15, Xie21}.
  \item Lusztig's based ring $J$ (also known as the asymptotic Hecke algebra) is defined with the assumptions that $(W,S,L)$ is bounded and that P1--P15 are valid; see \cite[Section~18]{Lus03}.
  \item In the case $L = \ell$, Xi \cite[Theorem~1.5]{Xi12} showed that the validity of Conjecture~\ref{conj-bound} implies the existence of the lowest two-sided cell.
      In that theorem, a description for this two-sided cell is also provided.
      Later, Xie \cite[Theorem 2.1]{Xie17} generalized this result to general $(W,S,L)$.
      In addition, in the case $L = \ell$, the existence of the lowest two-sided cell can be deduced from only the assumption that $(W,S,\ell)$ is bounded by an arbitrary number; see \cite{HHS}.
\end{enumerate}

The main result of this paper is the following:

\begin{Thm}\label{thm-fxyz}
For any positively weighted Coxeter group $(W, S, L)$ of finite rank, there exists a constant $N'_L(W)$ such that $N'_L(W)$ is a bound for $(W, S, L)$.
\end{Thm}

In general, the bound $N'_L(W)$ we found might not be optimal.
It is possible that $N'_L(W)$ is greater than the conjectural bound $N_L(W)$, even in the case $L = \ell$; see Example~\ref{ex-greater-bound}.
Nevertheless, as mentioned above, in many circumstances what really matters is the existence of a bound rather than its explicit value.
As a side note, in some special cases our approach is capable to give the conjectural bound. See Section~\ref{sec-eg} for more detials.

\subsection{Reduction of Theorem~\ref{thm-fxyz}} \label{subsec-overview-proof}

In this subsection we present a sketch of the proof of Theorem~\ref{thm-fxyz}, which will be reduced to Theorems~\ref{thm-exist-intersecting} and \ref{thm-intersectbound}.
We fix a positively weighted Coxeter group $(W,S,L)$ of finite rank.
The Bruhat order is denoted by $\le$.
We also write $w < w'$ if $w \le w'$ and $w \ne w'$.

Suppose $x, y \in W$ and $y = s_1 s_2 \cdots s_k$ $(s_i \in S)$ is a reduced expression.
For a sequence $I$ in $\{1,2, \dots, k\}$, say, $i_1 < i_2 < \cdots < i_p$, we write $z_I := x s_1 \cdots \widehat{s_{i_1}} \cdots \widehat{s_{i_p}} \cdots s_k$ (omitting the terms $s_{i_1}, s_{i_2}, \dots, s_{i_p}$) and
\[\xi_I := (v^{L(s_{i_1})} - v^{- L(s_{i_1})}) (v^{L(s_{i_2})} - v^{- L(s_{i_2})}) \cdots (v^{L(s_{i_p})} - v^{- L(s_{i_p})}).\]
By induction on $k$ one easily deduces that
\begin{equation}\label{eq-expansion-TxTy}
T_x T_y=\sum_I\xi_I T_{z_I},
\end{equation}
where $I$ ranges over sequences $i_1<i_2<\cdots<i_{p}$ in $\{1,2,\dots,k\}$ such that
\begin{equation}\label{eq-<}
xs_1 s_2\cdots \widehat{s_{i_1}}\cdots\widehat{s_{i_{n-1}}}\cdots s_{i_n-1} s_{i_n} < xs_1 s_2 \cdots \widehat{s_{i_1}}\cdots\widehat{s_{i_{n-1}}}\cdots s_{i_n -1}
\end{equation}
for all $n=1,2, \dots, p$ (the omitted terms on both sides of \eqref{eq-<} are $s_{i_1}, s_{i_2}, \dots, s_{i_{n-1}}$).

Let $L_m = \max \{L(s) \mid s \in S\}$.
As the weight function $L$ is positive on $S$, it is clear that
\[v^{- p L_m} \xi_I \in \mathbb{Z} [v^{-1}]\]
for any sequence $I$ of length $p$ appearing in \eqref{eq-expansion-TxTy}.
Therefore, to prove Theorem~\ref{thm-fxyz}, it suffices to show:
\begin{quotation}
  the lengths of all the sequences $I$ appearing in \eqref{eq-expansion-TxTy} for all $x,y\in W$ are uniformly bounded.
\end{quotation}
We divide the proof of this statement into two steps.
The first is to associate with each sequence $I$ an intersecting set (see Definition~\ref{def-intersecting}) of cardinality $p$ (Theorem~\ref{thm-exist-intersecting}), and the second is to show that the cardinalities of such sets are uniformly bounded (Theorem~\ref{thm-intersectbound}).

\begin{Thm}\label{thm-exist-intersecting}
Suppose $x, y \in W$ and $y = s_1 s_2 \cdots s_k$ is a reduced expression.
Suppose $1 \le i_1<i_2<\cdots<i_{p} \le k$ and \eqref{eq-<} holds for all $n=1,2,\dots,p$.
Then there exists an intersecting set $\mathfrak{Q}$  with $|\mathfrak{Q}|=p$.
\end{Thm}

Note that the statement of Theorem~\ref{thm-exist-intersecting}, as well as the definition of intersecting sets, is unrelated to the weight function $L$.

\begin{Thm}\label{thm-intersectbound}
For the fixed Coxeter group $(W,S)$ of finite rank, there is a constant $M$ such that $|\mathfrak{Q}|\le M$ for any intersecting set $\mathfrak{Q}$.
\end{Thm}

We remark that in the proof of Theorem~\ref{thm-intersectbound}, the famous Ramsey's theorem is applied to significantly simplify the problem.

By Theorem~\ref{thm-intersectbound}, the number $N'(W) :=\max\{|\mathfrak{Q}| \mid \mathfrak{Q} \text{ is intersecting}\}$ is finite.
Setting $N'_L(W) = N'(W)L_m$, Theorem~\ref{thm-fxyz} follows once we prove the two theorems above.

\subsection{Outline of the paper}
The rest of this paper is organized as follows.
In Section~\ref{sec-preliminaries}, we recollect some preliminary knowledge about the geometric representation $V$, the Tits cone, and the behavior of the hyperplanes inside the contragradient representation $V^*$.
In addition, we translate the relation \eqref{eq-<} into the geometric language of chambers and hyperplanes in $V^*$.
In Section~\ref{sec-constr-intersec}, we give the construction of the intersecting set associated with each sequence $I$ appearing in \eqref{eq-expansion-TxTy}.
In Section~\ref{sec-bound-intsect}, with the help of Ramsey's theorem, we prove that the cardinalities of intersecting sets are uniformly bounded.
In Section~\ref{sec-eg}, under the assumption $L=\ell$, we show that our upper bound coincides with the conjectural bound in the following cases: (1) finite; (2) affine; (3) $|S|=3$.
At the end of this section, we provide an example in the case $|S|=4$ showing that the bound we found might be bigger then the conjectural bound.
In the appendix, we include some interesting reformulations of some statements in this paper, provided by the anonymous referees.

\subsection*{Acknowledgments}
The authors would like to thank Professors Grant T.\ Barkley, Kei Yuen Chan, Junbin Dong, Tao Gui, Xuhua He, George Lusztig, Sian Nie, Nanhua Xi, Xun Xie for their helpful suggestions and enlightening discussions.
We also thank the anonymous referees for their valuable comments and additions.
The first author is grateful to Professors Jianpan Wang and Naihong Hu for their valuable advices.
The second author is supported by the National Natural Science Foundation of China (no.~12501037) and the Fundamental Research Funds for the Central Universities (no.~531118010972).

\section{Preliminaries}\label{sec-preliminaries}

\subsection{Hyperplanes, Tits cone, and Bruhat order} \label{subsec-Titscone}

We refer the readers to \cite{BB05, Bourbaki-Lie456, Hum, Kr} for this subsection.

Let $V$ be the $\mathbb{R}$-vector space with basis $\Delta=\{\alpha_s\mid s\in S\}$
endowed with a symmetric bilinear form $\operatorname{B}(-,-)$ defined by
\[\operatorname{B}(\alpha_s, \alpha_t) = \begin{cases}
    -\cos\frac{\pi}{m_{st}}, & \text{ if } m_{st} \ne \infty; \\
    -1, & \text{ if } m_{st} = \infty,
\end{cases}\]
where $m_{st}\in\{1,2,\cdots\}\cup\{\infty\}$ is the order of $st$.
The \textit{\textbf{geometric representation}} of $W$ on $V$ is given by
\begin{equation}\label{eq-ref}
sv = v-2 \operatorname{B}(\alpha_s,v)\alpha_s, \quad s \in S, \, v \in V.
\end{equation}
This representation is faithful and the $W$-action on $V$ preserves the bilinear form $\operatorname{B}$.

Let $\Phi=\{w \alpha_s \in V \mid w\in W, s\in S\}$ be the root system.
We have the partition $\Phi = \Phi^+ \sqcup \Phi^-$ where $\Phi^+$ is the set of positive roots (that is, the roots in $\Phi \cap \sum_{s \in S}\mathbb{R}_{\ge 0} \alpha_s$) and $\Phi^- := - \Phi^+$.
It is well known that the set $\Phi^+$ is in one-to-one correspondence with the set of \textit{\textbf {reflections}} $\bigcup_{w \in W} w S w^{-1}$.
For $\alpha \in \Phi^+$, say, $\alpha = w \alpha_s$ ($w \in W, \, s \in S$), we denote the corresponding reflection by
\[\sigma_\alpha := w s w^{-1}.\]
Then, for any $v \in V$, we have
\[\sigma_\alpha v = v - 2\operatorname{B} (v, \alpha) \alpha.\]
Thus, $\sigma_\alpha$ is independent of the choices of $w$ and $s$.

Let $V^*$ be the dual space of $V$, and $\langle -, -\rangle$ be the natural pairing $V^* \times V \to \mathbb{R}$.
The contragradient representation of $W$ on $V^*$ is characterized by
\[\langle wf, wv \rangle = \langle f, v\rangle, \quad w \in W, \ f \in V^*, \ v \in V.\]
For each root $\alpha \in \Phi$ we introduce the hyperplane
\(H_\alpha := \{ f \in V^* \mid \langle f, \alpha \rangle = 0 \}\).
Then, $H_\alpha = H_\beta$ if and only if $\alpha = \pm \beta$.
We set
\[\mathfrak{P} := \{H_\alpha \mid \alpha \in \Phi^+\}.\]
Since $w H_\alpha = H_{w \alpha}$, the $W$-action on $V^*$ induces a $W$-action on $\mathfrak{P}$.
For $H \in \mathfrak{P}$, we denote by $\alpha_H$ the positive root corresponding to $H$.
By abuse of notation, we also denote by $\sigma_H$ the reflection $\sigma_{\alpha_H}$.
Clearly, $\sigma_H$ fixes $H$ pointwise.

Let $C := \{f \in V^* \mid \langle f, \alpha_s \rangle > 0 \text{ for all } s \in S\}$ be the dominant chamber.
Its translations $wC$, $w \in W$, are called chambers.
All the chambers are disjoint to each other and they are in one-to-one correspondence with elements of $W$.
Each chamber $wC$ is an open simplicial cone with walls $wH_{\alpha_s}$, $s \in S$.

\begin{Rem}
    Since the chambers are simplicial cones, we often draw illustrative pictures on the unit sphere centered at the origin of $V^*$.
    On that sphere a chamber becomes a simplex (alcove).
    See \cite[Fig.~1--3 at the end of the paper]{Lus85} for such pictures for affine Weyl groups of types $\widetilde{A}_2$, $\widetilde{B}_2$ and $\widetilde{G}_2$.
\end{Rem}

Suppose $H \in \mathfrak{P}$ and $f, g \in V^*$.
If $\langle f,\alpha_H\rangle\langle g,\alpha_H\rangle<0$, that is, if $f$ and $g$ lie on different sides of $H$, then we say $H$ {\bfseries \itshape separates} $f$ and $g$.
Moreover, if $A,B \subset V^* \setminus H$ and all of the points in $A$ lie on one side of $H$ while those in $B$ lie on the other side, then we say $H$ {\bfseries \itshape separates} $A$ and $B$. The following result is well known.

\begin{Lem} \label{lem-chamber}
    Let $x, y \in W$, $s \in S$ be arbitrary.
    \begin{enumerate}
        \item The distance between $xC$ and $yC$, defined as the number of hyperplanes in $\mathfrak{P}$ separating $xC$ and $yC$, is equal to $\ell(y^{-1}x)$.
        \item The only hyperplane separating $xC$ and $xsC$ is $xH_{\alpha_s}$.
        \item The following are equivalent:
        \begin{enumerate}
            \item $xs < x$;
            \item $x \alpha_s \in \Phi^-$;
            \item the hyperplane $xH_{\alpha_s}$ separates $C$ and $xC$.
        \end{enumerate}
    \end{enumerate}
\end{Lem}

Let $\overline{C}$ be the closure of $C$ (with respect to the Euclidean topology).
The \textit{\textbf{Tits cone}} $U$ is defined by
\[U := \bigcup_{w \in W} w \overline{C}.\]
It is well known that $U$ and its interior $U^\circ$ are both convex cones in $V^*$ stable under the action of $W$.
Moreover, $\overline{C}$ is a fundamental domain for the $W$-action on $U$.

The closed cone $\overline{C}$ is a disjoint union of faces, namely, $\overline{C} = \sqcup_{I \subset S} C_I$, where
\[C_I = \{f \in V^* \mid \langle f, \alpha_s\rangle =0 \text{ for all } s \in I \text{ and } \langle f, \alpha_s \rangle > 0 \text{ for all } s \notin I\}.\]
If $f \in C_I$ for some $I \subset S$, then its stabilizer $\operatorname{Stab}_W(f)$ equals the standard parabolic subgroup $W_I$ generated by $I$.
In general, the closure of a chamber $w\overline{C}$ is partitioned into $w\overline{C} = \sqcup_{I \subset S} w C_I$, and the stabilizer of a point in $wC_I$ is the parabolic subgroup $wW_Iw^{-1}$.

\begin{Lem}\label{lem-stab}
    Let $f \in U$.
    The stabilizer $\operatorname{Stab}_W(f)$ is finite if and only if $f \in U^\circ$.
\end{Lem}

\begin{proof}
  See, for example, \cite[Proposition 2.2.5]{Kr}.
\end{proof}

In particular, since $\operatorname{Stab}_W(f)$ is trivial for any $f \in wC$, each chamber $w C$ is contained in $U^\circ$.

\begin{Cor} \label{cor-single-intersecting}
    For any $H \in \mathfrak{P}$, we have $H \cap U^\circ \ne \varnothing$.
\end{Cor}

\begin{proof}
    Suppose $H = H_{w \alpha_s} = w H_{\alpha_s}$, $w \in W$, $s \in S$.
    Then, $H \cap wC_{\{s\}} \ne \varnothing$.
    Let $f \in H \cap wC_{\{s\}}$.
    Then $\operatorname{Stab}_W(f) = \{e,wsw^{-1}\}$ is a finite subgroup of $W$.
    By Lemma~\ref{lem-stab}, we have $f \in U^\circ$.
\end{proof}

We will also need the following two lemmas which are essentially \cite[Propositions~4.5.4 and 4.5.5]{BB05}.

\begin{Lem} \label{lem-finite-dihedral}
   Suppose $P,Q \in \mathfrak{P}$, $P \ne Q$, then we have $P \cap Q \cap U^\circ \ne \varnothing$ if and only if $|\operatorname{B}(\alpha_P, \alpha_Q)|<1$.
\end{Lem}

\begin{proof}
It is a basic fact that the subgroup $W'$ generated by $\sigma_P, \sigma_Q$ is finite if and only if $|\operatorname{B}(\alpha_P, \alpha_Q)|<1$; see, for example, \cite[Proposition~4.5.4]{BB05}. It remains to prove that $P \cap Q \cap U^\circ \ne \varnothing$ if and only if $W'$ is finite.

If $P \cap Q \cap U^\circ \ne \varnothing$, then the reflections $\sigma_P$ and $\sigma_Q$ fix some point of $U^\circ$. Therefore, by Lemma~\ref{lem-stab}, $W'$ is a finite dihedral group.

Conversely, assume that $W'$ is finite.
We choose an element $f\in C\subset U^\circ$ and set $\overline{f}=\frac{1}{|W'|}\sum_{w\in W'}wf$.
Clearly, $f$ is fixed by $\sigma_P,\sigma_Q$ and hence $\overline{f}\in P\cap Q$.
As each $wf$ ($w \in W'$) lies in $U^\circ$ and $U^\circ$ is a convex cone, we have $\overline{f}\in U^\circ$.
Therefore, $P \cap Q \cap U^\circ \ne \varnothing$.
\end{proof}

\begin{Lem} \label{lem-cos-finite}
   If $|S|$ is finite, then the set
    \[\{c \in \mathbb{R} \mid |c|<1, \, c = \operatorname{B}(\alpha, \beta) \text{ for some } \alpha, \beta \in \Phi^+\}\]
    is finite.
\end{Lem}

\begin{proof}
    For any $\alpha \in \Phi^+$, there exists an element $w \in W$ such that $w\alpha=\alpha_s$ for some $s \in S$.
    Then, we have $\operatorname{B}(\alpha, \beta) = \operatorname{B}(\alpha_s, w\beta)$.
    Moreover, \cite[Proposition~4.5.5]{BB05} says that the set
    $$\{\operatorname{B}(\alpha_s,\gamma)\mid s\in S, \, \gamma\in\Phi, \, |\operatorname{B}(\alpha_s,\gamma)|<1\}$$
    is finite if $|S|$ is finite. Therefore, the set in question is finite.
\end{proof}

\subsection{Intersecting sets}

For any subsets $A, B$ of $V^*$, we define
\[A\dc B := A\cap B\cap U^\circ.\]

\begin{Def} \leavevmode
    \begin{enumerate}
        \item Let $\mathfrak{P}_1, \mathfrak{P}_2 \subset \mathfrak{P}$ be two subsets.
           We say that $\mathfrak{P}_1$ and $\mathfrak{P}_2$ {\bf \itshape intersect pairwise} (or equivalently, $\mathfrak{P}_1$ intersects pairwise with $\mathfrak{P}_2$) if $H_1 \dc H_2 \ne \varnothing$ for any $H_1 \in \mathfrak{P}_1$ and $H_2 \in \mathfrak{P}_2$.
        \item If moreover, $\mathfrak{P}_1 = \{P\}$ consists of a single hyperplane, we also say $P$ \emph{intersects pairwise} with $\mathfrak{P}_2$, or equivalently, $\mathfrak{P}_2$ intersects pairwise with $P$.
    \end{enumerate}
\end{Def}

By convention, we understand that $\varnothing$ intersects pairwise with any subset of $\mathfrak{P}$.

\begin{Def}\label{def-intersecting}
A subset $\mathfrak{Q}$ of $\mathfrak{P}$ is called {\bf \itshape intersecting} if $\mathfrak{Q}$ intersect pairwise with itself, that is, if $H_1\dc H_2\ne \varnothing$ for any $H_1,H_2\in\mathfrak{Q}$. In particular, $\varnothing$ is intersecting.
\end{Def}

\begin{Rem} \label{rmk-I(1)} \leavevmode
  \begin{enumerate}
      \item A single hyperplane $H \in \mathfrak{P}$ is regarded to be intersecting because $H \cap U^\circ \ne \varnothing$ (see Corollary~\ref{cor-single-intersecting}).
      \item If $\mathfrak{Q}$ is intersecting then $w(\mathfrak{Q}) := \{wP \mid P \in \mathfrak{Q}\}$ is also intersecting for any $w \in W$.
      This is because $U^\circ$ is preserved by the $W$-action.
  \end{enumerate}

\end{Rem}

For $f,g \in V^*$, we denote the (closed) segment between $f$ and $g$ by
\[[f,g] := \{\lambda f + (1- \lambda)g \mid 0 \le \lambda \le 1\}.\]
The following lemma is obvious in affine geometry.

\begin{Lem}\label{lem-0ptthm}
Let $A$ be a convex subset of $V^*$ and  $P\in\mathfrak{P}$. If $\langle f_1, \alpha_P \rangle \langle f_2, \alpha_P \rangle \le0$ for some $f_1,f_2\in A$ (in particular, if $P$ separates $f_1$ and $f_2$), then $P\cap A\ne \varnothing$.
\end{Lem}

\begin{proof}
The intermediate value theorem implies that $[f_1,f_2]\cap P\ne \varnothing$ (in fact, the set $[f_1,f_2]\cap P$ consists of a single point).
Meanwhile, we have $[f_1,f_2]\subset A$ by the convexity of $A$.
Thus, $P \cap A \ne \varnothing$.
\end{proof}

The following two results are easy consequences of Lemma~\ref{lem-0ptthm}.

\begin{Lem}\label{lem-isolateintersect}
Let $f_1,f_2\in U^\circ$ and $P,P'\in\mathfrak{P}$ such that $f_1, f_2 \notin P \cup P'$. Assume that
\begin{enumerate}
    \item $P'$ separates $f_1$ and $f_2$ while $P$ does not;
    \item there exists a point $g \in P' \cap U^\circ$ such that $P$ separates $g$ and $\{f_1, f_2\}$.
\end{enumerate}
Then, $P\dc P'\ne \varnothing$.
\end{Lem}

See Figure~\ref{fig-lem-isolateintersect} for an illustration.

\begin{figure}[ht]
    \centering
    \begin{tikzpicture}
      \draw (-2,0.2) -- (2,0.2);
      \draw (0,-1) -- (0,2);
      \coordinate [label = right:$P$] (P) at (2,0.2);
      \coordinate [label = right:$P'$] (Pp) at (0,2);
      \coordinate [label=above:$f_1$] (f1) at (-1,0.8);
      \coordinate [label=above:$f_2$] (f2) at (1,0.8);
      \fill (0,-0.5) circle (2pt) node [left] {$g$};
    \end{tikzpicture}
    \caption{Illustration for Lemma~\ref{lem-isolateintersect}} \label{fig-lem-isolateintersect}
\end{figure}

\begin{proof}
By Lemma~\ref{lem-0ptthm} and the fact that $U^\circ$ is convex, we have $[f_1, f_2] \dc P' \ne \varnothing$.
Moreover, $P$ separates $g$ and the point $[f_1, f_2] \dc P'$.
Since $g$ and $[f_1, f_2] \dc P'$ are located in $P'$, it is clear that $P \dc P' \ne \varnothing$ by Lemma~\ref{lem-0ptthm} again.
\end{proof}

\begin{Lem} \label{lem-determine-intersect}
    Let $P, P' \in \mathfrak{P}$, $f_1, f_2, g_1, g_2 \in U^\circ$ such that $f_i,g_i \notin P \cup P'$ ($i = 1, 2$).
    Suppose
    \begin{enumerate}
        \item $P'$ separates $f_1$ and $f_2$ while $P$ does not;
        \item $P'$ separates $g_1$ and $g_2$ while $P$ does not;
        \item $P$ separates $\{f_1, f_2\}$ and $\{g_1, g_2\}$.
    \end{enumerate}
    Then, $P \dc P' \ne \varnothing$.
\end{Lem}

See Figure~\ref{fig-determine-intersect} for an illustration.

\begin{figure}[ht]
    \centering
    \begin{tikzpicture}
        \draw (-2,0) -- (2,0);
        \draw (0,-1.5) -- (0,1.5);
        \coordinate [label = right:$P$] (P) at (2,0);
        \coordinate [label = right:$P'$] (P') at (0,1.5);
        \coordinate [label = center:$f_1$] (f1) at (-1,0.7);
        \coordinate [label = center:$g_1$] (g1) at (-1,-0.7);
        \coordinate [label = center:$f_2$] (f2) at (1,0.7);
        \coordinate [label = center:$g_2$] (g2) at (1,-0.7);

        \node  (or) at (3.7,0) {or};

        \draw (5,0) -- (9,0);
        \draw (7,-1.5) -- (7,1.5);
        \coordinate [label = right:$P$] (oP) at (9,0);
        \coordinate [label = right:$P'$] (oP') at (7,1.5);
        \coordinate [label = center:$f_1$] (of1) at (6,0.7);
        \coordinate [label = center:$g_2$] (og2) at (6,-0.7);
        \coordinate [label = center:$g_1$] (og1) at (8,-0.7);
        \coordinate [label = center:$f_2$] (of2) at (8,0.7);
    \end{tikzpicture}
    \caption{Illustration for Lemma~\ref{lem-determine-intersect}. Roughly speaking, the assumptions (1), (2), and (3) are saying that the four points $f_1, f_2, g_1, g_2$ are located in the four quadrants of $V^* \setminus (P \cup P')$ respectively.}
    \label{fig-determine-intersect}
\end{figure}

\begin{proof}
By Lemma~\ref{lem-0ptthm} and the fact that $U^\circ$ is convex, we have $[g_1, g_2] \dc P' \ne \varnothing$.
Moreover, $P$ separates $[g_1, g_2] \dc P'$ and $\{f_1, f_2\}$.
Now the assertion follows from Lemma~\ref{lem-isolateintersect}.
\end{proof}

When applying the above two lemmas, we need to determine whether two points are separated by a certain hyperplane.
The following lemma is useful for this purpose.
Recall the notation $\sigma_H$ indicates the reflection with respect to a hyperplane $H$, and $\alpha_H$ indicates the corresponding positive root.

\begin{Lem} \label{lem-y-sigmay-sameside}
    Let $P, H\in \mathfrak{P}$, $f, g \in U^\circ$ such that $f,g \notin P \cup H$.
    Suppose
    \begin{enumerate}
        \item $f$ and $g$ lie on the same side of $H$;
        \item $f$ and $g$ lie on different sides of $P$;
        \item $f$ and $\sigma_H f$ lie on different sides of $P$.
    \end{enumerate}
    Then, $g$ and $\sigma_Hg$ lie on the same side of $P$.
\end{Lem}

See Figure~\ref{fig-y-sigmay-sameside} for an illustration.

\begin{figure}[ht]
    \centering
    \begin{tikzpicture}
        \draw (-2,0) -- (2.5,0);
        \draw (60:2) -- (240:2);
        \coordinate [label=right:$H$] (H) at (2.5,0);
        \coordinate [label=right:$P$] (P) at (60:2);
        \node (f) at (0,1.3) {$f$};
        \node (g) at (1.5,0.8) {$g$};
        \node (sf) at (0,-1.3) {$\sigma_Hf$};
        \node (sg) at (1.5,-0.8) {$\sigma_Hg$};
    \end{tikzpicture}
    \caption{Illustration for Lemma~\ref{lem-y-sigmay-sameside}}
    \label{fig-y-sigmay-sameside}
\end{figure}

\begin{proof}
    Replacing $\alpha_H$ (or $\alpha_P$) by its negative $-\alpha_H$ (or $-\alpha_P$) when necessary, we may assume that $\langle f, \alpha_H \rangle > 0$ and $\langle f, \alpha_P\rangle >0$.
    Note that
    \begin{align*}
        \langle \sigma_H f, \alpha_P\rangle & = \langle f, \sigma_H \alpha_P \rangle \\
        & = \langle f, \alpha_P - 2 \operatorname{B} (\alpha_P, \alpha_H) \alpha_H \rangle \\
        & = \langle f, \alpha_P\rangle - 2 \operatorname{B} (\alpha_P, \alpha_H) \langle f, \alpha_H\rangle.
    \end{align*}
    By the assumption (3), we have $\langle \sigma_H f, \alpha_P \rangle < 0$.
    It yields that $\operatorname{B} (\alpha_P, \alpha_H) > 0$.
    Similarly, we also have
    \[\langle \sigma_Hg, \alpha_P \rangle = \langle g, \alpha_P \rangle - 2 \operatorname{B} (\alpha_P, \alpha_H) \langle g, \alpha_H\rangle.\]
    By the assumptions (1) and (2), we have $\langle g, \alpha_H\rangle > 0$ and $\langle g, \alpha_P \rangle < 0$.
    Therefore, $\langle \sigma_Hg, \alpha_P \rangle < 0$.
    In conclusion, $g$ and $\sigma_Hg$ lie on the same side of $P$.
\end{proof}

\subsection{The path from \texorpdfstring{$xC$}{xC} to \texorpdfstring{$xyC$}{xyC}} \label{subsec-path-notations}

As in Section~\ref{subsec-overview-proof}, we fix $x, y \in W$ and a reduced expression $y = s_1 \cdots s_k$.
Suppose $1 \le i_1 < \cdots < i_p \le k$ such that \eqref{eq-<} holds for each $n=1,2,\dots,p$.
In the proof of Theorem~\ref{thm-exist-intersecting} (Section~\ref{sec-constr-intersec}), we adopt the following notations ($1 \le n\le p$):
\begin{align*}
    x_n & :=x s_1 s_2 \cdots s_{i_n-1}, \\
    \sigma_n & := x_ns_{i_n}x_n^{-1}, \text{ the reflection sending $x_n C$ to $x_n s_{i_n}C$}, \\
    H_n &\in\mathfrak{P}  \text{ :  the hyperplane corresponding to the reflection $\sigma_n$}, \\
    e_n &:= \sigma_{n}\cdots\sigma_{1}.
\end{align*}
Moreover, we regard $e_0 = e$.

Consider the sequence of chambers
\[xC, \quad xs_1C, \quad xs_1s_2C, \quad \dots, \quad xs_1 \cdots s_k C.\]
For convenience, we temporarily denote $\widetilde{x}_j := xs_1 \cdots s_{j}$ ($0 \le j \le k$, regard $\widetilde{x}_0 = x$).
By Lemma~\ref{lem-chamber}, two adjacent chambers $\widetilde{x}_{j-1}C$ and $\widetilde{x}_{j}C$ are separated by a unique hyperplane $\widetilde{x}_{j-1} H_{\alpha_{s_{j}}}$.
Moreover, since $s_1 \cdots s_k$ is a reduced expression, the distance between $\widetilde{x}_{j'}C$ and $\widetilde{x}_jC$ ($j' < j$) is exactly $j - j'$.
One easily deduces by induction that the hyperplanes separating $\widetilde{x}_{j'}C$ and $\widetilde{x}_jC$ are exactly
\[\widetilde{x}_{j'} H_{\alpha_{s_{j'+1}}}, \quad \widetilde{x}_{j'+1} H_{\alpha_{s_{j'+2}}}, \quad \dots, \quad \widetilde{x}_{j-1} H_{\alpha_{s_{j}}}.\]
It follows immediately that the chambers $xC, \ \widetilde{x}_1C, \ \dots, \ \widetilde{x}_{j-1}C$ lie on one side of the hyperplane $\widetilde{x}_{j-1} H_{\alpha_{s_j}}$ while the others $\widetilde{x}_{j}C, \ \widetilde{x}_{j+1}C, \ \dots, \ xyC$ lie on the other side.
(The path from $xC$ to $xyC$ in the ``alcove picture" is illustrated in Figure~\ref{fig-path}.)
In particular, we have the following fact.

\begin{Lem}\label{lem-x-xt-Ht}
For each $1\le n\le p$, the chambers $x_1C$ and $x_nC$ lie on the same side of $H_n$, while $x_n s_{i_n}C$ lies on the other side.
\end{Lem}

\begin{proof}
    In our notation, $x_1 = \widetilde{x}_{i_1 - 1}$, $x_n = \widetilde{x}_{i_n -1}$, $x_ns_{i_n} = \widetilde{x}_{i_n}$, $H_n = \widetilde{x}_{i_n -1} H_{\alpha_{s_{i_n}}}$.
    The lemma follows from the above discussion.
\end{proof}

\begin{figure}[ht]
\begin{tikzpicture}
\path [draw] (0,1.5) -- (0,-1.5);
\path [draw] (-3.06066, 1.06066) -- (-0.93934, -1.06066);
\path [draw] (0.93934, -1.06066) -- (3.06066, 1.06066);
\draw ($(-5,0) + (90:0.7) + (120:1)$) -- ($(-5,0) + (90:0.7) + (120:-2.3)$);
\draw ($(5,0) + (150:0.7) + (120:1)$) -- ($(5,0) + (150:0.7) + (120:-2)$);
\coordinate [label=right:$xH_{\alpha_{s_1}}$] (Hs1) at (-5.5,1.5);
\coordinate [label=right:$xyH_{\alpha_{s_k}}$] (lastplane) at (3.8,1.4);
\fill [pattern=dots, pattern color=lightgray] ($(-5,0) + (90:0.7)$) -- ($(-5,0) + (210:0.7)$) -- ($(-5,0) + (330:0.7)$);
\fill [pattern=dots, pattern color=lightgray] ($(5,0) + (270:0.7)$) -- ($(5,0) + (30:0.7)$) -- ($(5,0) + (150:0.7)$);
\fill [pattern=dots, pattern color=lightgray] (0,0.6) -- ($(0,0.6) + (210:1.2)$) -- (0,-0.6);
\fill [pattern=dots, pattern color=lightgray] ($(-2,0) + (135:0.6)$) -- ($(-2,0) + (135:0.6) + (255:1.2)$) -- ($(-2,0) + (315:0.6)$);
\fill [pattern=dots, pattern color=lightgray] ($(2,0) + (45:0.6)$) -- ($(2,0) + (45:0.6) + (165:1.2)$) -- ($(2,0) + (225:0.6)$);
\coordinate [label=above:$H_1$] (h1) at (-3.06066, 1.06066);
\coordinate [label=above:$H_n$] (ht) at (0,1.5);
\coordinate [label=above:$H_p$] (hn) at (3.06066, 1.06066);
\coordinate [label=center:$xC$] (xc) at (-5,0);
\coordinate [label=center:$xyC$] (xyc) at (5,0);
\coordinate [label=center:$\cdots$] (cdots1) at (-1.2,1);
\coordinate [label=center:$\cdots$] (cdots2) at (1.2,1);
\coordinate [label=center:$\cdots$] (d1) at (-3.5,0.5);
\coordinate [label=center:$\cdots$] (d2) at (3.5,0.5);
\coordinate [label=below left:$x_1C$] (u1c) at (-1.8,0);
\coordinate [label=left:$x_nC$] (utc) at (0.1,0);
\coordinate [label=above left:$x_pC$] (unc) at (2.28,-0.1);
\draw[->](2,0) -- +(-45:0.5cm);
\draw[->](0,0) -- +(0:0.5cm);
\draw[->](-2,0) -- +(45:0.5cm);
\draw [->] ($(-5,0) + (30:0.35)$) -- + (30:0.5cm);
\coordinate [label=above:$s_1$] (s1) at ($(-4.9,0.1) + (30:0.35)$);
\coordinate [label=above:$s_{i_1}$] (si1) at (-2,0.1);
\coordinate [label=above right:$s_{i_n}$] (st) at (0,0);
\coordinate [label=right:$s_{i_p}$] (sn) at (2,0);
\end{tikzpicture}
\caption{The path $xC = \widetilde{x}_0C\rightarrow \widetilde{x}_1C\rightarrow \widetilde{x}_2C\rightarrow\cdots\rightarrow \widetilde{x}_{k}C=xyC$}\label{fig-path}
\end{figure}

The following lemma translates \eqref{eq-<} into a statement on the relationship of $H_n$, $x_1C$, and $e_{n-1}C$.

\begin{Lem}\label{lem-H-sep-x-e}
For each $1\le n\le p$, the hyperplane $H_n$ separates $x_1C$ and $e_{n-1}C$.
\end{Lem}

\begin{proof}
For each $1\le n\le p$, let $x_n':=xs_1\cdots \widehat{s_{i_1}}\cdots\widehat{s_{i_{n-1}}}\cdots \widehat{s_{i_n}}$ (in particular, $x_1' = x_1$).
Then \eqref{eq-<} becomes $x_n' s_{i_n} < x_n'$.
By Lemma~\ref{lem-chamber}, this is equivalent to saying that the hyperplane $x_n' H_{\alpha_{s_{i_n}}}$ separates $C$ and $x_n' C$.
This is transformed into the following statement by applying the action of $e_{n-1}$:
\[\text{the hyperplane $e_{n-1}x_n' H_{\alpha_{s_{i_n}}}$ separates $e_{n-1}C$ and $e_{n-1}x_n' C$.}\]
In addition, we have $e_{n-1} x_n' = \sigma_{n-1} \cdots \sigma_1x_n' = x_n$ by an easy computation.
Note also that $x_n H_{\alpha_{s_{i_n}}} = H_n$ in our notation.
Therefore, the statement above reads:
\[\text{the hyperplane $H_n$ separates $e_{n-1}C$ and $x_n C$.}\]
In view of Lemma~\ref{lem-x-xt-Ht}, $x_1C$ and $x_nC$ lie on the same side of $H_n$.
Hence, $H_n$ separates $e_{n-1}C$ and $x_1 C$.
\end{proof}

As the reflection $\sigma_n$ swaps $e_{n-1}C$ and $e_{n}C$, we have the following corollary.

\begin{Cor} \label{cor-H-sep-x-e}
     The two chambers $x_1C$ and $e_nC$ lie on the same side of $H_n$.
\end{Cor}

\begin{proof}
    It is clear that $H_n$ separates $e_{n-1}C$ and $e_{n}C$.
    Then the assertion follows immediately from Lemma~\ref{lem-H-sep-x-e}.
\end{proof}

\section{Proof of Theorem~\ref{thm-exist-intersecting}} \label{sec-constr-intersec}

We fix $x, y \in W$, a reduced expression $y = s_1 \cdots s_k$, and a sequence $i_1, \dots, i_p$ such that \eqref{eq-<} holds for each $n = 1, \dots, p$.
We retain the notations $x_n,  \sigma_n,  e_n,  H_n$ as in Section~\ref{subsec-path-notations}.

Recall that Theorem~\ref{thm-exist-intersecting} claims the existence of an intersecting set $\mathfrak{Q}$ of cardinality $p$.
Rather than constructing directly such a set, we will construct inductively a sequence of intersecting sets $\mathfrak{Q}_n$ ($1 \le n \le p$) of cardinality $n$ satisfying stronger properties (Theorem~\ref{thm-exist-intersecting-2}) and set $\mathfrak{Q}=\mathfrak{Q}_p$.
In the statement of this theorem, the first two properties ${\bf C}(n)$ and ${\bf I}(n)$ are what we really need, while the others ensure that our inductive construction can proceed.

\begin{Thm} \label{thm-exist-intersecting-2}
For each $n = 1, \dots, p$, there exists a subset $\mathfrak{Q}_n \subset \mathfrak{P}$ such that the following are satisfied:
    \begin{enumerate}
        \item [${\bf C}(n)$:] $|\mathfrak{Q}_n|=n$;
        \item [${\bf I}(n)$:] $\mathfrak{Q}_n$ is intersecting;
        \item [${\bf H}(n)$:] $H_n \in \mathfrak{Q}_n$;
        \item [${\bf S}(n)$:] $P$ separates $x_1C$ and $x_ns_{i_n}C$ for any $P \in \mathfrak{Q}_n$;
        \item [${\bf U}(n)$:] $P$ does not separate $x_1C$ and $e_nC$ for any $P \in \mathfrak{Q}_n$.
    \end{enumerate}
    In particular, Theorem~\ref{thm-exist-intersecting} holds.
\end{Thm}

\begin{Rem}
    For convenience of memorization, the letters ${\bf C},{\bf I},{\bf H},{\bf S},{\bf U}$ stand for cardinality, intersecting, the letter H in $H_n$, separate, unilateral, respectively.
\end{Rem}

\subsection{The construction of intersecting sets \texorpdfstring{$\mathfrak{Q}_n$}{Qn}}\label{construct-Qn}

Our strategy to construct $\mathfrak{Q}_n$ is an induction on $n$.
For $n = 1$, the requirements ${\bf C}(1)$ and ${\bf H}(1)$ force $\mathfrak{Q}_1 = \{H_1\}$.
The property ${\bf I}(1)$ holds automatically (see Remark~\ref{rmk-I(1)}(1)), while ${\bf S}(1)$ and ${\bf U}(1)$ follow from Lemma~\ref{lem-x-xt-Ht} and Corollary~\ref{cor-H-sep-x-e}, respectively.

Suppose that $n \ge 2$ and that we have $\mathfrak{Q}_{n-1}$ in hand satisfying ${\bf C}(n-1)$, ${\bf I}(n-1)$, ${\bf H}(n-1)$, ${\bf S}(n-1)$ and ${\bf U}(n-1)$.
We will partition $\mathfrak{Q}_{n-1}$ into a disjoint union of two subsets, say, $\mathfrak{Q}_{n-1} = \mathfrak{A} \sqcup \mathfrak{B}$, and define
\begin{equation}\label{eq-Qn}
    \mathfrak{Q}_n = \sigma_{n}(\mathfrak{A}) \cup \mathfrak{B} \cup \{H_n\}.
\end{equation}

Now consider a hyperplane $P \in \mathfrak{Q}_{n-1}$ and suppose $P$ separates $x_1C$ and $e_nC$.
Then, $P$ does not satisfy the property {\bf U}($n$).
Thus, if we hope the $\mathfrak{Q}_n$ constructed in \eqref{eq-Qn} to satisfy ${\bf U}(n)$, then the hyperplane $P$ can not belong to $\mathfrak{B}$.
In other words, $\mathfrak{A}$ must contain the following set,
\[\mathfrak{A}_0 := \{P \in \mathfrak{Q}_{n-1} \mid P \text{ separates } x_1C \text{ and } e_nC\}.\]
Unfortunately, the set $\sigma_{n}(\mathfrak{A}_0) \cup (\mathfrak{Q}_{n-1} \setminus \mathfrak{A}_0) \cup \{H_n\}$ is not adequate to be the desired $\mathfrak{Q}_n$.
This is because $\sigma_{n}(\mathfrak{A}_0)$ might not intersect pairwise with $\mathfrak{Q}_{n-1} \setminus \mathfrak{A}_0$ and this violates the requirement ${\bf I}(n)$.
Therefore, $\mathfrak{A}$ must also contain all the hyperplanes in $\mathfrak{Q}_{n-1} \setminus \mathfrak{A}_0$ which are not intersecting pairwise with $\sigma_{n}(\mathfrak{A}_0)$.
For this, we define inductively
\[\mathfrak{A}_{i+1} := \{P \in \mathfrak{Q}_{n-1} \setminus (\mathfrak{A}_0 \sqcup  \dots \sqcup \mathfrak{A}_{i}) \mid P \text{ does not intersect pairwise with } \sigma_{n} \mathfrak{A}_{i}\}.\]
As mentioned, the requirement ${\bf I}(n)$, together with the containment $\mathfrak{A}_0 \subset \mathfrak{A}$, forces $\mathfrak{A}_1 \subset \mathfrak{A}$.
By the same reason, ${\bf I}(n)$ and the containment $\mathfrak{A}_1 \subset \mathfrak{A}$ force $\mathfrak{A}_2 \subset \mathfrak{A}$.
Inductively, $\mathfrak{A}_i \subset \mathfrak{A}$ for each $i$.

Note that we regard the empty set intersects pairwise with any hyperplane.
Therefore, if $\mathfrak{A}_l = \varnothing$ for some $l \ge 0$, then $\mathfrak{A}_i = \varnothing$ for any $i \ge l$.
As $\mathfrak{Q}_{n-1}$ is a finite set, such $l$ exists.
At last we define
\begin{equation} \label{eq-Qn'}
    \begin{split}
        \mathfrak{A} &:= \bigsqcup_{i=0}^\infty \mathfrak{A}_i = \mathfrak{A}_0 \sqcup \mathfrak{A}_1 \sqcup \dots \sqcup \mathfrak{A}_{l-1}, \\
        \mathfrak{B} &:= \mathfrak{Q}_{n-1} \setminus \mathfrak{A}, \\
        \mathfrak{Q}_n &:= \sigma_{n}(\mathfrak{A}) \cup \mathfrak{B} \cup \{H_n\}.
    \end{split}
\end{equation}
See Figure~\ref{fig-construction} for an illustration.

\begin{figure}[ht]
    \centering
    \begin{tikzpicture}
        \draw [thick] (-1.1,0) -- (3,0);
        \coordinate [label=above:$\mathfrak{B}_\sigma$] (Bsigma) at (2.2,0);
        \coordinate [label=above:$\mathfrak{A}_{l-1}$] (Al) at (-0.4,0);
        \draw (1.4,0.1) -- (1.4,-0.1);
        \draw (-0.8,0.1) -- (-0.8,-0.1);
        \draw [thick] (3,0.2) -- (3,-0.2);
        \draw [thick] (0,0.2) -- (0,-0.2);
        \draw [thick] (-6.5,0.2) -- (-6.5,-0.2);
        \draw [thick] (-6.5,0) -- (-1.9,0);
        \coordinate [label=above:$\mathfrak{A}_0$] (A0) at (-5.45,0);
        \coordinate [label=above:$\mathfrak{A}_{1}$] (A1) at (-3.8,0);
        \coordinate [label=above:$\mathfrak{A}_2$] (A2) at (-2.7,0);
        \node (d) at (-1.5,0) {$\dots$};
        \node (Qn-1) at (-7.5,0) {$\mathfrak{Q}_{n-1}$:};
        \draw (-4.4,0.1) -- (-4.4,-0.1);
        \draw (-3.2,0.1) -- (-3.2,-0.1);
        \draw (-2.2,0.1) -- (-2.2,-0.1);
        \draw [decorate, decoration={brace, amplitude=5pt, mirror}, yshift=-3pt] (-6.45,-0.3) -- (-0.05,-0.3) node [midway, below=5pt] {$\mathfrak{A}$};
        \draw [decorate, decoration={brace, amplitude=5pt, mirror}, yshift=-3pt] (0.05,-0.3) -- (2.95,-0.3) node [midway, below=5pt] {$\mathfrak{B}$};
    \end{tikzpicture}
    \caption{Illustration for the partition $\mathfrak{Q}_{n-1} = \mathfrak{A} \sqcup \mathfrak{B}$. The set $\mathfrak{B}_\sigma$ will be defined in Section~\ref{subsec-3.2}, and the fact that $\mathfrak{B}_\sigma \subset \mathfrak{B}$ will be given in Proposition~\ref{prop-Bsigma-in-B}.}
    \label{fig-construction}
\end{figure}

By definition, we have the property ${\bf H}(n)$, that is, $H_n \in \mathfrak{Q}_n$.
In the rest of this section we will see that the set $\mathfrak{Q}_n$ satisfies ${\bf C}(n)$, ${\bf I}(n)$, ${\bf S}(n)$, ${\bf U}(n)$, which completes the proof of Theorem~\ref{thm-exist-intersecting-2}.

\subsection{Some preparations} \label{subsec-3.2}

Before we prove the properties ${\bf C}(n)$, ${\bf I}(n)$, ${\bf S}(n)$, ${\bf U}(n)$, we need some preliminary results.

\begin{Lem} \label{lem-Hn-notin-Qn-1} \leavevmode
    \begin{enumerate}
        \item We have $H_n \notin \mathfrak{Q}_{n-1}$.
        \item For any $P  \in \mathfrak{Q}_{n-1}$, $P$ does not separate $x_nC$  and $x_ns_{i_{n}}C$.
        \item For any $P  \in \mathfrak{Q}_{n-1}$, $P$ separates $x_1C$ and $x_nC$.
    \end{enumerate}
\end{Lem}

\begin{proof}
    We have seen in Lemma~\ref{lem-H-sep-x-e} that $H_n$ separates $x_1 C$ and $e_{n-1}C$. This violates ${\bf U}(n-1)$ and  hence $H_n \notin \mathfrak{Q}_{n-1}$. This proves (1).

    Note that $H_n$ is the unique hyperplane separating $x_nC$ and $x_ns_{i_n}C$.
    Thus, (2) follows from (1).

    For (3), notice that $P$ separates $x_1C$ and $x_{n-1}s_{i_{n-1}}C$ by ${\bf S}(n-1)$.
    By the discussion before Lemma~\ref{lem-x-xt-Ht}, a hyperplane separating $x_1C$ and $x_{n-1}s_{i_{n-1}}C$ must separate $x_1C$ and $x_n C$ as well.
\end{proof}

\begin{Lem} \label{lem-Qn-1-int-Hn}
    The hyperplane $H_n$ intersects pairwise with $\mathfrak{Q}_{n-1}$.
\end{Lem}

\begin{proof}
    Suppose $P \in \mathfrak{Q}_{n-1}$.
    By Lemma~\ref{lem-H-sep-x-e} and ${\bf U}(n-1)$, we have
    \begin{enumerate}
        \item [(a)] $H_n$ separates $x_1C$ and $e_{n-1}C$ while $P$ does not.
    \end{enumerate}
    By Lemma~\ref{lem-chamber} and Lemma~\ref{lem-Hn-notin-Qn-1}(2), we have
    \begin{enumerate}
        \item [(b)] $H_n$ separates $x_nC$ and $x_ns_{i_{n}}C$ while $P$ does not.
    \end{enumerate}
    It also holds that (see Lemma~\ref{lem-Hn-notin-Qn-1}(3))
    \begin{enumerate}
        \item [(c)] $P$ separates $x_1C$ and $x_{n}C$.
    \end{enumerate}
    By Lemma~\ref{lem-determine-intersect}, we deduce from (a), (b), and (c) that $P \dc H_n \ne \varnothing$.
\end{proof}

\begin{Lem} \label{lem-A0}
    Suppose $P \in \mathfrak{A}_0$. Then,
    \begin{enumerate}
        \item $x_1C, \ \sigma_nx_1C, \ e_{n-1}C$ lie on the same side of $P$, while $x_nC, \ x_ns_{i_n}C, \ e_nC$ lie on the other side;
        \item $x_1C, \ \sigma_nx_1C, \ e_{n}C$ lie on the same side of $\sigma_n P$, while $x_nC, \ x_ns_{i_n}C, \ e_{n-1}C$ lie on the other side.
    \end{enumerate}
\end{Lem}

See Figure~\ref{fig-A0} for an illustration.

\begin{figure}[ht]
    \centering
    \begin{tikzpicture}
        \draw (-3,0) -- (3,0);
        \draw (60:3) -- (240:3);
        \draw (120:3) -- (300:3);
        \coordinate [label=right:$H_n$] (Hn) at (3,0);
        \coordinate [label=right:$P \in \mathfrak{A}_0$] (P) at (300:3);
        \coordinate [label=right:$\sigma_nP$] (sP) at (60:3);

        \fill [pattern=dots, pattern color=lightgray] ($(150:2) + (90:0.6)$) -- ($(150:2) + (210:0.6)$) -- ($(150:2) + (330:0.6)$) -- cycle;
        \fill [pattern=dots, pattern color=lightgray] ($(270:2) + (90:0.6)$) -- ($(270:2) + (210:0.6)$) -- ($(270:2) + (330:0.6)$) -- cycle;
        \fill [pattern=dots, pattern color=lightgray] ($(2,0.3) + (90:0.6)$) -- ($(2,0.3) + (210:0.6)$) -- ($(2,0.3) + (330:0.6)$) -- cycle;
        \fill [pattern=dots, pattern color=lightgray] ($(90:2) + (270:0.6)$) -- ($(90:2) + (150:0.6)$) -- ($(90:2) + (30:0.6)$) -- cycle;
        \fill [pattern=dots, pattern color=lightgray] ($(210:2) + (270:0.6)$) -- ($(210:2) + (150:0.6)$) -- ($(210:2) + (30:0.6)$) -- cycle;
        \fill [pattern=dots, pattern color=lightgray] ($(2,-0.3) + (270:0.6)$) -- ($(2,-0.3) + (150:0.6)$) -- ($(2,-0.3) + (30:0.6)$) -- cycle;

        \node (x1) at (150:2) {$x_1C$};
        \node (en) at (90:2) {$e_nC$};
        \node (sx1) at (210:2) {$\sigma_nx_1C$};
        \node (en-1) at (270:2) {$e_{n-1}C$};
        \node (xn) at (2,0.3) {$x_nC$};
        \node (sxn) at (2,-0.3) {$x_ns_{i_n}C$};
    \end{tikzpicture}
    \caption{Illustration for Lemma~\ref{lem-A0}}
    \label{fig-A0}
\end{figure}

\begin{proof}
    By Corollary~\ref{cor-H-sep-x-e}, we have
    \begin{enumerate}
        \item [(a)] $e_nC$ and $x_1C$ lie on the same side of $H_n$.
    \end{enumerate}
    Since $P \in \mathfrak{A}_0$, we have
    \begin{enumerate}
        \item [(b)] $e_nC$ and $x_1C$ lie on different sides of $P$.
    \end{enumerate}
    By ${\bf U}(n-1)$, $e_{n-1}C$ and $x_1C$ lie on the same side of $P$. This, as well as the fact (b), yields that
    \begin{enumerate}
        \item [(c)] $e_nC$ and $e_{n-1}C$ ($=\sigma_n e_nC$) lie on different sides of $P$.
    \end{enumerate}
    By Lemma~\ref{lem-y-sigmay-sameside}, the facts (a), (b) and (c) imply that
    \begin{enumerate}
        \item [(d)] $x_1C$ and $\sigma_n x_1C$ lie on the same side of $P$.
    \end{enumerate}
    Combining the facts (b), (d), ${\bf U}(n-1)$, ${\bf S}(n-1)$, and Lemma~\ref{lem-Hn-notin-Qn-1}(2), the point (1) is proved.
    Applying the reflection $\sigma_n$, we obtain the point (2).
\end{proof}

We define
\[\mathfrak{B}_\sigma := \{P \in \mathfrak{Q}_{n-1} \mid \sigma_{n}P \text{ separates } x_1C \text{ and } e_nC\}.\]
As we will see, $\mathfrak{B}_\sigma$ plays a role dual to that of $\mathfrak{A}_0$.

\begin{Lem} \label{lem-A0-cap-Bsigma-empty}
    $\mathfrak{A}_0 \cap \mathfrak{B}_\sigma = \varnothing$.
\end{Lem}

\begin{proof}
    Suppose $P \in \mathfrak{A}_0$.
    Then, by Lemma~\ref{lem-A0}(2), $\sigma_n P$ does not separate $x_1C$ and $e_nC$.
    This contradicts the definition of $\mathfrak{B}_\sigma$.
    Therefore, $\mathfrak{A}_0 \cap \mathfrak{B}_\sigma = \varnothing$.
\end{proof}

\begin{Lem} \label{lem-Bsigma}
    Suppose $P \in \mathfrak{B}_\sigma$. Then,
    \begin{enumerate}
        \item $x_1C, \ e_{n-1}C, \ e_nC$ lie on the same side of $P$, while $\sigma_nx_1C,\ x_nC,\ x_ns_{i_n}C$ lie on the other side;
        \item $x_1C, \ x_nC, \ x_ns_{i_n}C$ lie on the same side of $\sigma_n P$, while $\sigma_nx_1C, \ e_{n-1}C, \ e_{n}C$ lie on the other side.
    \end{enumerate}
\end{Lem}

See Figure~\ref{fig-Bsigma} for an illustration.

\begin{figure}[ht]
    \centering
    \begin{tikzpicture}
        \draw (-3,0) -- (3,0);
        \draw (60:3) -- (240:3);
        \draw (120:3) -- (300:3);
        \coordinate [label=right:$H_n$] (Hn) at (3,0);
        \coordinate [label=right:$P \in \mathfrak{B}_\sigma$] (P) at (60:3);
        \coordinate [label=right:$\sigma_nP$] (sP) at (300:3);

        \fill [pattern=dots, pattern color=lightgray] ($(150:2) + (90:0.6)$) -- ($(150:2) + (210:0.6)$) -- ($(150:2) + (330:0.6)$) -- cycle;
        \fill [pattern=dots, pattern color=lightgray] ($(270:2) + (90:0.6)$) -- ($(270:2) + (210:0.6)$) -- ($(270:2) + (330:0.6)$) -- cycle;
        \fill [pattern=dots, pattern color=lightgray] ($(2,0.3) + (90:0.6)$) -- ($(2,0.3) + (210:0.6)$) -- ($(2,0.3) + (330:0.6)$) -- cycle;
        \fill [pattern=dots, pattern color=lightgray] ($(90:2) + (270:0.6)$) -- ($(90:2) + (150:0.6)$) -- ($(90:2) + (30:0.6)$) -- cycle;
        \fill [pattern=dots, pattern color=lightgray] ($(210:2) + (270:0.6)$) -- ($(210:2) + (150:0.6)$) -- ($(210:2) + (30:0.6)$) -- cycle;
        \fill [pattern=dots, pattern color=lightgray] ($(2,-0.3) + (270:0.6)$) -- ($(2,-0.3) + (150:0.6)$) -- ($(2,-0.3) + (30:0.6)$) -- cycle;

        \node (x1) at (150:2) {$e_nC$};
        \node (en) at (90:2) {$x_1C$};
        \node (sx1) at (210:2) {$e_{n-1}C$};
        \node (en-1) at (270:2) {$\sigma_nx_1C$};
        \node (xn) at (2,0.3) {$x_nC$};
        \node (sxn) at (2,-0.3) {$x_ns_{i_n}C$};
    \end{tikzpicture}
    \caption{Illustration for Lemma~\ref{lem-Bsigma}}
    \label{fig-Bsigma}
\end{figure}

\begin{proof}
    Since $P \in \mathfrak{B}_\sigma$ and $\mathfrak{A}_0 \cap \mathfrak{B}_\sigma = \varnothing$ (Lemma~\ref{lem-A0-cap-Bsigma-empty}), we have $P \notin \mathfrak{A}_0$ and hence
    \begin{enumerate}
        \item [(a)] $x_1C$ and $e_nC$ lie on the same side of $P$.
    \end{enumerate}
    By the definition of $\mathfrak{B}_\sigma$, $\sigma_{n}P$ separates $x_1C$ and $e_nC$, which is equivalent to
    \begin{enumerate}
        \item [(b)] $P$ separates $\sigma_nx_1C$ and $e_{n-1}C$.
    \end{enumerate}
    Combining the facts (a), (b), ${\bf U}(n-1)$, and Lemma~\ref{lem-Hn-notin-Qn-1}(2)(3), the point (1) is proved.
    Applying the reflection $\sigma_n$, we obtain the point (2).
\end{proof}

\subsection{Proof of \texorpdfstring{${\bf C}(n)$, ${\bf I}(n)$, ${\bf S}(n)$, and ${\bf U}(n)$}{C(n), I(n), S(n), and U(n)}}

Recall that we assume the existence of $\mathfrak{Q}_{n-1}$ and define $\mathfrak{Q}_n$ as in \eqref{eq-Qn'}.
In this subsection we prove the properties ${\bf C}(n)$, ${\bf I}(n)$, ${\bf S}(n)$ and ${\bf U}(n)$.

\subsubsection{Proof of \texorpdfstring{${\bf C}(n)$}{C(n)}} \label{subsubsec-Cn}

\begin{Lem} \label{lem-sigmaAi-cap-B-empty}
    For each $i \ge 0$, it holds that $\sigma_n (\mathfrak{A}_i) \cap \mathfrak{Q}_{n-1} = \varnothing$.
    In particular, $\sigma_n (\mathfrak{A}_i) \cap \mathfrak{B} = \varnothing$.
\end{Lem}

\begin{proof}
    Suppose first $P \in\mathfrak{A}_0$.
    If $\sigma_nP \in \mathfrak{Q}_{n-1}$, then $\sigma_nP$ does not separate $x_1C$ and $e_{n-1}C$ by ${\bf U}(n-1)$.
    But this contradicts Lemma~\ref{lem-A0}(2).
    Hence, $\sigma_n (\mathfrak{A}_0) \cap \mathfrak{Q}_{n-1} = \varnothing$.

    Suppose now $P \in \mathfrak{A}_i$ ($i \ge 1$) and $\sigma_nP \in \mathfrak{Q}_{n-1}$.
    By the definition of $\mathfrak{A}_i$, there is a hyperplane $Q \in \mathfrak{A}_{i-1}$ such that $\sigma_nP \dc Q =\varnothing$.
    This contradicts ${\bf I}(n-1)$.
\end{proof}

We have seen that $H_n \notin \mathfrak{Q}_{n-1}$ in Lemma~\ref{lem-Hn-notin-Qn-1}(1).
Hence, $H_n \notin \sigma_n(\mathfrak{A}) \cup \mathfrak{B}$.
In view of Lemma~\ref{lem-sigmaAi-cap-B-empty}, the union
\[\sigma_{n}(\mathfrak{A}) \cup \mathfrak{B} \cup \{H_n\}\]
is a disjoint union.
Therefore, $|\mathfrak{Q}_n| = |\mathfrak{A}| + |\mathfrak{B}| + 1 = |\mathfrak{Q}_{n-1}| + 1 = n$.
The property ${\bf C}(n)$ is proved.

\subsubsection{Proof of \texorpdfstring{${\bf I}(n)$}{I(n)}}

\begin{Lem} \label{lem-sigmaA-B-int}
    The two sets of hyperplanes $\sigma_n(\mathfrak{A})$ and $\mathfrak{B}$ intersect pairwise.
\end{Lem}

\begin{proof}
    Suppose $A \in \mathfrak{A}$, $B \in \mathfrak{B}$.
    We may assume that $A \in \mathfrak{A}_i$.
    If $\sigma_nA \dc B =\varnothing$, then $B \in \mathfrak{A}_{i+1}$ by the definition of $\mathfrak{A}_{i+1}$, which contradicts $B \in \mathfrak{B}$.
    Therefore $\sigma_n(\mathfrak{A})$ and $\mathfrak{B}$ intersect pairwise.
\end{proof}

\begin{Lem} \label{cor-Qn-1-int-Hn}
    The hyperplane $H_n$ intersects pairwise with $\sigma_{n}(\mathfrak{A}) \cup \mathfrak{B}$.
\end{Lem}

\begin{proof}
    If $P \in \mathfrak{B}$, then $P \dc H_n \ne \varnothing$ by Lemma~\ref{lem-Qn-1-int-Hn}.
    If $P = \sigma_nP'$ for some $P' \in \mathfrak{A}$, then $P \dc H_n=\sigma_nP' \dc H_n = \sigma_n(P' \dc H_n) \ne \varnothing$ by Lemma~\ref{lem-Qn-1-int-Hn} again.
\end{proof}

At last, $\sigma_n(\mathfrak{A})$ is intersecting since $\mathfrak{A}$ is intersecting (by ${\bf I}(n-1)$) and $\sigma_n$ preserves the property of being intersecting (see Remark~\ref{rmk-I(1)}(2)).
The set $\mathfrak{B}$ is also intersecting by ${\bf I}(n-1)$.
These, together with the results in Lemma~\ref{lem-sigmaA-B-int} and Lemma~\ref{cor-Qn-1-int-Hn}, show that $\mathfrak{Q}_n$ ($= \sigma_{n}(\mathfrak{A}) \cup \mathfrak{B} \cup \{H_n\}$) is intersecting.
The property ${\bf I}(n)$ is proved.

\subsubsection{Proof of \texorpdfstring{${\bf S}(n)$ and ${\bf U}(n)$}{S(n) and U(n)}}

To prove ${\bf S}(n)$ and ${\bf U}(n)$, we need the following result.

\begin{Prop}\label{prop-Bsigma-in-B}
    We have $\mathfrak{B}_\sigma \subset \mathfrak{B}$.
    In other words, $\mathfrak{A} \cap \mathfrak{B}_\sigma = \varnothing$.
\end{Prop}

The proof of Proposition~\ref{prop-Bsigma-in-B} is postponed until Section~\ref{subset-pf-Bsigma-in-B} as it is much more involved than all the above lemmas.
Assuming its validity, we can now prove ${\bf S}(n)$ and ${\bf U}(n)$.

\begin{proof}[Proof of ${\bf S}(n)$]
For any $P \in \mathfrak{Q}_n$, we need to show that $P$ separates $x_1C$ and $x_ns_{i_{n}}C$.
There are three cases.

\emph{Case 1.} If $P = H_n$, it has been done in Lemma~\ref{lem-x-xt-Ht}.

\emph{Case 2.} If $P \in \mathfrak{B}$, then in particular $P \in \mathfrak{Q}_{n-1}$.
By Lemma~\ref{lem-Hn-notin-Qn-1}(2)(3), $P$ separates $x_1C$ and $x_{n}s_{i_{n}}C$.

\emph{Case 3.} Suppose $P \in \sigma_n(\mathfrak{A})$, that is, $\sigma_nP \in \mathfrak{A}$.
Then, by Proposition~\ref{prop-Bsigma-in-B}, we have $\sigma_nP \notin \mathfrak{B}_\sigma$.
By the definition of $\mathfrak{B}_\sigma$, it holds that
\begin{enumerate}
    \item [(a)] $x_1C$ and $e_nC$ lie on the same side of $P$.
\end{enumerate}
By ${\bf U}(n-1)$, we see that $x_1C$ and $e_{n-1}C$ lie on the same side of $\sigma_nP$.
Applying the reflection $\sigma_n$, we obtain
\begin{enumerate}
    \item [(b)] $\sigma_n x_1C$ and $e_{n}C$ lie on the same side of $P$.
\end{enumerate}
Therefore, by (a) and (b), we have
\begin{enumerate}
    \item [(c)] $x_1C$ and $\sigma_n x_1C$ lie on the same side of $P$.
\end{enumerate}
By Lemma~\ref{lem-Hn-notin-Qn-1}(3), $\sigma_nP$ separates $x_1C$ and $x_{n}C$.
Therefore, by applying the reflection $\sigma_n$, we have
\begin{enumerate}
    \item [(d)] $P$ separates $\sigma_n x_1C$ and $x_n s_{i_n}C$.
\end{enumerate}
The facts (c) and (d) imply that $P$ separates $x_1C$ and $x_n s_{i_n}C$.

The property ${\bf S}(n)$ is proved.
\end{proof}

\begin{proof}[Proof of ${\bf U}(n)$]
For any $P \in \mathfrak{Q}_n$, we need to show that $x_1C$ and $e_nC$ lie on the same side of $P$.
There are three cases.

\emph{Case 1.} If $P = H_n$, it has been done in Corollary~\ref{cor-H-sep-x-e}.

\emph{Case 2.} Suppose $P \in \mathfrak{B}$.
Then $P \notin \mathfrak{A}_0$.
By the definition of $\mathfrak{A}_0$, $P$ does not separate $x_1C$ and $e_nC$.

\emph{Case 3.} Suppose $P \in \sigma_n(\mathfrak{A})$. Then $\sigma_nP\in\mathfrak{A}$, and hence $\sigma_n P\not\in\mathfrak{B}_\sigma$ by Proposition~\ref{prop-Bsigma-in-B}.
Thus, by the definition of $\mathfrak{B}_\sigma$, $P$ does not separate $x_1C$ and $e_nC$.

The property ${\bf U}(n)$ is proved.
\end{proof}

We have proved Theorem~\ref{thm-exist-intersecting-2} modulo Proposition~\ref{prop-Bsigma-in-B}.

\subsection{Proof of Proposition~\ref{prop-Bsigma-in-B}} \label{subset-pf-Bsigma-in-B}

To prove Proposition~\ref{prop-Bsigma-in-B}, we need to show that:
\begin{equation}\label{mainstep}
    \mathfrak{A}_i\cap\mathfrak{B}_\sigma = \varnothing \text{ for any } i \ge 0.
\end{equation}
The case $i = 0$ has been settled in Lemma~\ref{lem-A0-cap-Bsigma-empty}.

Suppose now $i \ge 1$ and $P_i$ is an arbitrary hyperplane in $\mathfrak{A}_i$.
By the definition of $\mathfrak{A}_j$ $(0<j\le i)$, there exist hyperplanes $P_j\in\mathfrak{A}_j$ $(0\le j < i)$ such that
\begin{equation} \label{eq-assumption-Pj}
 \sigma_n P_{j}\dc P_{j+1} = \varnothing \text{ for any }  j = 0,1, \dots, i-1.
\end{equation}
In the rest of this subsection we fix the index $i$ and a choice of the above sequence $P_0, P_1, \dots, P_i$.

Note that these hyperplanes are intersecting by ${\bf I}(n-1)$.
In addition, $P_j \dc H_n \ne \varnothing$ for any $j$ by Lemma~\ref{lem-Qn-1-int-Hn}.
We will need the following lemmas concerning the relative position of these intersections.

\begin{Lem} \label{lem-Ai-same-side}
    Suppose $i \ge 2$.
    For any $j = 2, \dots, i$, the two (nonempty) sets $P_{j-2}\dc H_n$ and $P_{j}\dc H_n$ lie on the same side of $\sigma_nP_{j-1}$.
\end{Lem}

\begin{proof}
    Let $f_1 \in P_{j}\dc H_n$, $f_2 \in P_{j-2}\dc H_n$ and $g \in P_{j-2} \dc P_{j}$ be arbitrary.
    By the condition \eqref{eq-assumption-Pj}, it holds that $f_1, f_2, g \notin \sigma_n P_{j-1}$.
    Suppose that $\sigma_nP_{j-1}$ separates $f_1$ and $f_2$.
    Then, $g$ is separated from either $f_1$ or $f_2$ by $\sigma_nP_{j-1}$.
    If $\sigma_nP_{j-1}$ separates $g$ and $f_1$, then $P_j \dc \sigma_n P_{j-1} \ne \varnothing$ by Lemma~\ref{lem-0ptthm} as both $g$ and $f_1$ are points in $P_j \cap U^\circ$ (see Figure~\ref{fig-Ai-same-side} for an illustration).
    This contradicts \eqref{eq-assumption-Pj}.
    If otherwise $\sigma_nP_{j-1}$ separates $g$ and $f_2$, then for the same reason we have $P_{j-2} \dc \sigma_n P_{j-1} \ne \varnothing$ which is also absurd.
    To conclude, $f_1$ and $f_2$ lie on the same side of $\sigma_nP_{j-1}$.
\begin{figure}[ht]
    \centering
    \begin{tikzpicture}
      \path [draw, name path = Hline] (-4,0) -- (2.5,0);
      \coordinate [label=right:$H_n$] (h) at (2.5,0);
      \path [draw, dashed, name path = hPline] (-1.7,-0.9) -- (0,3);
      \coordinate [label = right:$\sigma_n P_{j-1}$] (hp) at (-1.7,-0.8);
      \path [draw, name path = H1line] (-3, -0.9) -- (-1,3);
      \coordinate [label = left:$P_{j-2}$] (h1) at (-2.9,-0.8);
      \path [draw, name path = H2line] (-2, 3) -- (1,-0.9);
      \coordinate [label = right:$P_{j}$] (h2) at (0.9,-0.8);
      \fill [name intersections = {of = H1line and Hline, by = {[label=135:$f_2$]f1}}] (f1) circle (2pt);
      \fill [name intersections = {of = H2line and Hline, by = {[label=45:$f_1$]f2}}] (f2) circle (2pt);
      \fill [name intersections = {of = H2line and H1line, by = {[label=left:$g$]f2}}] (f2) circle (2pt);
    \end{tikzpicture}
    \caption{Illustration for the proof of Lemma~\ref{lem-Ai-same-side}} \label{fig-Ai-same-side}
\end{figure}
\end{proof}

For a hyperplane $P \in \mathfrak{P}$, we denote by $P^{x_1,-}$ the open half space of $V^* \setminus P$ not containing $x_1C$, that is,
\[P^{x_1,-} := \{f \in V^* \mid P \text{ separates } f \text{ and } x_1C \}.\]
Let $\alpha_P \in \Phi^+$ be the root corresponding to $P$ and $f_0 \in x_1C$ be arbitrary.
Then $P^{x_1,-}$ is characterized by
\[P^{x_1,-} = \{f \in V^* \mid \langle f, \alpha_P \rangle \langle f_0, \alpha_P \rangle < 0\}.\]

Let $\alpha_0, \alpha_1, \dots, \alpha_i$ and $\alpha_H \in \Phi^+$ be the roots corresponding to $P_0, P_1, \dots, P_i$ and $H_n$, respectively.
Then, $\sigma_n P_j$ ($j = 0,1, \dots, i$) is characterized by
\[\sigma_n P_j = \{f \in V^* \mid \langle f, \sigma_n \alpha_j \rangle = 0\}.\]

\begin{Lem} \label{lem-affine-geometry}
    Let $j$ be fixed.
    Suppose there exist $g_1, g_2 \in H_n^{x_1, -} \cap (\sigma_n P_j)^{x_1, -}$ such that $P_j$ separates $g_1$ and $g_2$.
    Then $P_j \cap H_n^{x_1, -} =P_j \cap (\sigma_n P_j)^{x_1, -}$.
\end{Lem}

See Figure~\ref{fig-affine-geometry} for an illustration.

\begin{figure}[ht]
    \centering
    \begin{tikzpicture}
        \draw (-2.5,0) -- (2.5,0);
        \draw (60:3) -- (240:3);
        \draw [dashed] (120:3) -- (300:3);
        \coordinate [label=above:$H_n$] (Hn) at (2.5,0);
        \coordinate [label=right:$P_j$] (P) at (300:3);
        \coordinate [label=right:$\sigma_nP_j$] (sP) at (60:3);

        \fill [pattern=dots, pattern color=lightgray] ($(150:2) + (90:0.6)$) -- ($(150:2) + (210:0.6)$) -- ($(150:2) + (330:0.6)$) -- cycle;

        \node (x1) at (150:2) {$f_0 \in x_1C$};
        \node (en-1) at (270:2) {$g_1$};
        \node (sxn) at (330:2) {$g_2$};
        \fill ($(270:2)+(30:1)$) circle (2pt);
        \coordinate [label=right:$g$] (g) at ($(270:2)+(30:1)$);

        \node (or) at (3.7,0) {or};

        \draw (4.5,0) -- (9.5,0);
        \draw ($(60:3) + (7,0)$) -- ($(240:3)+ (7,0)$);
        \draw [dashed] ($(120:3) + (7,0)$) -- ($(300:3) + (7,0)$);
        \coordinate [label=above:$H_n$] (Hn') at (9.5,0);
        \coordinate [label=right:$P_j$] (P') at ($(300:3)+ (7,0)$);
        \coordinate [label=right:$\sigma_nP_j$] (sP') at ($(60:3)+ (7,0)$);

        \fill [pattern=dots, pattern color=lightgray] ($(90:2) + (270:0.6)+ (7,0)$) -- ($(90:2) + (150:0.6)+ (7,0)$) -- ($(90:2) + (30:0.6)+ (7,0)$) -- cycle;

        \node (x1') at ($(90:2)+ (7,0)$) {$f_0 \in x_1C$};
        \node (en-1') at ($(270:2)+ (7,0)$) {$g_1$};
        \node (sxn') at ($(330:2)+ (7,0)$) {$g_2$};
        \fill ($(270:2)+(30:1)+ (7,0)$) circle (2pt);
        \coordinate [label=right:$g$] (g') at ($(270:2)+(30:1)+ (7,0)$);
    \end{tikzpicture}
    \caption{Illustration for Lemma~\ref{lem-affine-geometry}}
    \label{fig-affine-geometry}
\end{figure}

\begin{proof}
    Let $\Lambda=H_n^{x_1, -} \cap (\sigma_n P_j)^{x_1, -}$. Then $\Lambda$ is convex. Since  $g_1, g_2 \in \Lambda$ and $P_j$ separates $g_1$ and $g_2$, Lemma~\ref{lem-0ptthm} implies that there is a point $g\in P_j\cap\Lambda$ (in fact, $\{g\} = [g_1, g_2] \cap P_j$), and hence
    \begin{equation} \label{eq-affine-geometry-1}
         \langle g, \alpha_H \rangle \langle f_0, \alpha_H \rangle < 0 \quad \text{and} \quad \langle g, \sigma_n\alpha_{j} \rangle \langle f_0, \sigma_n\alpha_{j} \rangle < 0
    \end{equation}
    for any $f_0 \in x_1C$.
    In particular, $g \notin H_n\cap\sigma_n P_j$.
    Note that $H_n \ne P_j$ by Lemma~\ref{lem-Hn-notin-Qn-1}(1).
    Therefore, $H_n \cap \sigma_nP_j$ is a linear subspace of $V^*$ of codimension $2$ and hence a subspace of $P_j$ of codimension $1$.
    Thus, $P_j = (H_n \cap \sigma_nP_j) \oplus \mathbb{R} g$.
    Now we have
    \begin{equation} \label{eq-affine-geometry-3}
        P_j \cap H_n^{x_1, -} = \{f \in (H_n \cap \sigma_nP_j) \oplus \mathbb{R} g \mid \langle f, \alpha_H \rangle \langle f_0, \alpha_H \rangle < 0\}.
    \end{equation}
    By the first inequality in \eqref{eq-affine-geometry-1} and the fact that $\langle f', \alpha_H \rangle = 0$ for any $f' \in H_n$, it is clear from \eqref{eq-affine-geometry-3} that
    \[P_j \cap H_n^{x_1, -} = (H_n \cap \sigma_nP_j) + \mathbb{R}_{>0} g.\]
    Similarly, from the second inequality in \eqref{eq-affine-geometry-1} and the fact that $\langle f' , \sigma_n \alpha_j \rangle = 0$ for any $f' \in \sigma_n P_j$, we also deduce that
    \[P_j \cap (\sigma_n P_j)^{x_1, -} = (H_n \cap \sigma_nP_j) + \mathbb{R}_{>0} g.\]
    Therefore, $P_j \cap H_n^{x_1, -} = P_j \cap (\sigma_n P_j)^{x_1, -}.$
\end{proof}

\begin{Lem} \label{lem-Ai-other-side}
    For any $j =1, \dots, i$, we have $P_j \dc P_{j-1} \subset H_n^{x_1, -}$.
\end{Lem}

\begin{proof}
    We proceed with induction on $j$.
    Suppose first that $j = 1$.
    Note that $x_1C$ and $e_{n-1}C$ lie on the same side of $P_1$ (by ${\bf U}(n-1)$), while $x_1C$ and $x_nC$ are separated by $P_1$ (by Lemma~\ref{lem-Hn-notin-Qn-1}(3)).
    Therefore,
    \begin{enumerate}
        \item [(a)] $P_1$ separates $e_{n-1}C$ and $x_nC$.
    \end{enumerate}
    In view of Lemma~\ref{lem-A0}(2), it holds that
    \begin{enumerate}
        \item [(b)] $e_{n-1}C$ and $x_nC$ are contained in the convex set $(\sigma_n P_0)^{x_1,-}\cap U^\circ$.
    \end{enumerate}
    Therefore, by (a), (b) and Lemma~\ref{lem-0ptthm},  we have $P_1\dc(\sigma_n P_0)^{x_1,-}\ne \varnothing$.
    It follows that $P_1 \cap U^\circ \subset (\sigma_n P_0)^{x_1,-}$ for, otherwise, there exist two points in $P_1 \cap U^\circ$ separated by $\sigma_n P_0$ and this forces $P_1 \dc \sigma_n P_0 \ne \varnothing$ which contradicts \eqref{eq-assumption-Pj}.
    In particular,
    \begin{enumerate}
        \item [(c)] $P_1 \dc P_0 \subset (\sigma_n P_0)^{x_1,-} \cap P_0$.
    \end{enumerate}
    Note also that we have seen the following facts in Lemmas~\ref{lem-x-xt-Ht}, \ref{lem-H-sep-x-e}, and \ref{lem-A0}:
    \begin{enumerate}
        \item [(d)] $P_0$ separates $e_{n-1}C$ and $x_n s_{i_n} C$;
        \item [(e)] $e_{n-1}C, \ x_n s_{i_n} C \subset H_n^{x_1,-} \cap (\sigma_n P_0)^{x_1,-}$.
    \end{enumerate}
    By Lemma~\ref{lem-affine-geometry}, we deduce from (d) and (e) that $(\sigma_n P_0)^{x_1,-} \cap P_0 = H_n^{x_1,-} \cap P_0$.
    From this and the fact (c), we obtain $P_1 \dc P_0 \subset H_n^{x_1,-}$ as desired (see Figure~\ref{fig-Ai-other-side-1} for an illustration).

\begin{figure}[ht]
    \centering
    \begin{tikzpicture}
        \draw (-3,0) -- (5,0);
        \draw [dashed] (60:3) -- (240:3.5);
        \path [draw, name path = P0line] (120:3) -- (300:3.5);
        \path [draw, name path = P1line] (0.3, -3) -- (3.4,2.5);
        \coordinate [label=right:$P_1 \in \mathfrak{A}_1$] (P1) at (3.4,2.5);
        \coordinate [label=right:$H_n$] (Hn) at (5,0);
        \coordinate [label=right:$P_0 \in \mathfrak{A}_0$] (P) at (300:3.5);
        \coordinate [label=right:$\sigma_nP_0$] (sP) at (60:3);

        \fill [name intersections = {of = P0line and P1line, by = {[label=right:$P_1 \dc P_0$]g}}] (g) circle (2pt);

        \fill [pattern=dots, pattern color=lightgray] ($(150:2) + (90:0.6)$) -- ($(150:2) + (210:0.6)$) -- ($(150:2) + (330:0.6)$) -- cycle;
        \fill [pattern=dots, pattern color=lightgray] ($(270:2) + (90:0.6)$) -- ($(270:2) + (210:0.6)$) -- ($(270:2) + (330:0.6)$) -- cycle;
        \fill [pattern=dots, pattern color=lightgray] ($(4,0.3) + (90:0.6)$) -- ($(4,0.3) + (210:0.6)$) -- ($(4,0.3) + (330:0.6)$) -- cycle;
        \fill [pattern=dots, pattern color=lightgray] ($(4,-0.3) + (270:0.6)$) -- ($(4,-0.3) + (150:0.6)$) -- ($(4,-0.3) + (30:0.6)$) -- cycle;

        \node (x1) at (150:2) {$x_1C$};
        \node (en-1) at (270:2) {$e_{n-1}C$};
        \node (xn) at (4,0.3) {$x_nC$};
        \node (sxn) at (4,-0.3) {$x_ns_{i_n}C$};
    \end{tikzpicture}
    \caption{Illustration for the proof of Lemma~\ref{lem-Ai-other-side}, case $j=1$}
    \label{fig-Ai-other-side-1}
\end{figure}

    Now we consider the case $j \ge 2$.
    We take arbitrarily
    \[g_- \in P_{j-1} \dc P_{j-2}, \quad g_+ \in P_j \dc P_{j-1}, \quad f_1 \in P_j \dc H_n, \quad f_2 \in P_{j-2} \dc H_n.\]
    By induction hypothesis, we have $g_- \in H_n^{x_1,-}$.
    Suppose $g_+ \notin H_n^{x_1,-}$.
    Then,
    \[\langle g_+, \alpha_H \rangle \langle g_-, \alpha_H \rangle \le 0.\]
    Since $\sigma_n \alpha_{j-1} = \alpha_{j-1} -2 \operatorname{B}(\alpha_{j-1}, \alpha_H) \alpha_H$ and $g_+, g_-\in P_{j-1}$,
    we have
    \begin{equation}\label{ineq1}
        \langle g_+, \sigma_n \alpha_{j-1} \rangle \langle g_-, \sigma_n \alpha_{j-1} \rangle = 4 \operatorname{B}(\alpha_{j-1}, \alpha_H)^2 \langle g_+, \alpha_{H} \rangle \langle g_-, \alpha_{H} \rangle \le 0.
    \end{equation}
    Note that $g_+$ and $f_1$ lie on the same side of $\sigma_n P_{j-1}$ because $g_+, f_1 \in P_j$ and $P_j \dc \sigma_n P_{j-1} = \varnothing$ (see \eqref{eq-assumption-Pj}), while $f_1$ and $f_2$ also lie on the same side of $\sigma_n P_{j-1}$ by Lemma~\ref{lem-Ai-same-side}.
    Hence, $g_+$ and $f_2$ lie on the same side of $\sigma_n P_{j-1}$, that is,
    \begin{equation}\label{ineq2}
        \langle g_+, \sigma_n \alpha_{j-1} \rangle \langle f_2, \sigma_n \alpha_{j-1} \rangle > 0.
    \end{equation}
Combining \eqref{ineq1} and \eqref{ineq2} yields
    \[ \langle f_2, \sigma_n \alpha_{j-1} \rangle \langle g_-, \sigma_n \alpha_{j-1} \rangle \le 0.\]
    Since $f_2, g_- \in P_{j-2} \cap U^\circ$, this implies $P_{j-2} \dc \sigma_n P_{j-1} \ne \varnothing$ by Lemma~\ref{lem-0ptthm}, which contradicts \eqref{eq-assumption-Pj}.
    Therefore, $g_+ \in H_n^{x_1, -}$ and hence $P_j \dc P_{j-1} \subset H_n^{x_1,-}$ as desired (see Figure~\ref{fig-Ai-other-side-2} for an illustration).
\begin{figure}[ht]
    \centering
    \begin{tikzpicture}
      \path [draw, name path = Hline] (-3,0) -- (6,0);
      \coordinate [label=right:$H_n$] (h) at (6,0);
      \path [draw, name path = Pline] (-1,-3) -- (1,3);
      \path [draw, name path = H1line] (0,3) -- (3,-2.5);
      \path [draw, name path = H2line] (-1.5,-2.7) -- (5.2,0.7);
      \draw [dashed] (-1,3) -- (1,-3);
      \coordinate [label=right:$P_{j-2}$] (h2) at (5.2,0.7);
      \coordinate [label=right:$P_j$] (h1) at (3,-2.5);
      \coordinate [label=right:$P_{j-1}$] (p) at (1,3);
      \coordinate [label=left:$\sigma_n P_{j-1}$] (sigmahp) at (-1,3);
      \fill [name intersections = {of = H1line and Hline, by = f1}] (f1) circle (2pt);
      \coordinate [label=45:$f_1$] (f1label) at (f1);
      \fill [name intersections = {of = H2line and Hline, by = f2}] (f2) circle (2pt);
      \coordinate [label=135:$f_2$] (f2label) at (f2);
      \fill [name intersections = {of = H1line and Pline, by = g1}] (g1) circle (2pt);
      \coordinate [label=right:$g_+$] (g1label) at (g1);
      \fill [name intersections = {of = H2line and Pline, by = g2}] (g2) circle (2pt);
      \coordinate [label=135:$g_-$] (g2label) at (g2);

      \fill [pattern=dots, pattern color=lightgray] ($(150:2) + (90:0.6)$) -- ($(150:2) + (210:0.6)$) -- ($(150:2) + (330:0.6)$) -- cycle;
      \node (x1) at (150:2) {$x_1C$};
    \end{tikzpicture}
    \caption{Illustration for the proof of Lemma~\ref{lem-Ai-other-side}, case $j \ge 2$}
    \label{fig-Ai-other-side-2}
\end{figure}
\end{proof}

With the above preparations in hand, we are able to complete the proof of Proposition~\ref{prop-Bsigma-in-B}.
\begin{proof}[Proof of Proposition~\ref{prop-Bsigma-in-B}]
Recall that $P_i$ ($i \ge 1$) is an arbitrary hyperplane in $\mathfrak{A}_i$ and $P_{i-1} \in \mathfrak{A}_{i-1}$ is chosen such that $P_{i-1} \dc \sigma_n P_i = \varnothing$.
Suppose further that $P_i \in \mathfrak{B}_\sigma$ (see Figure~\ref{fig-Pi-Bsigma} for an illustration).
Then, by Lemma~\ref{lem-H-sep-x-e} and Lemma~\ref{lem-Bsigma},
\begin{enumerate}
    \item [(a)] $P_i$ separates $e_{n-1}C$ and $\sigma_n x_1C$;
    \item [(b)] $e_{n-1}C, \ \sigma_n x_1C \subset H_{n}^{x_1, -} \cap (\sigma_n P_i)^{x_1, -}$.
\end{enumerate}
By Lemma~\ref{lem-affine-geometry}, the two facts above imply $P_i \cap H_n^{x_1,-} = P_i \cap (\sigma_n P_i)^{x_1,-}$.
Therefore, in view of Lemma~\ref{lem-Ai-other-side}, we have
\[P_i \dc P_{i-1} \subset P_i \cap H_n^{x_1,-} = P_i \cap (\sigma_n P_i)^{x_1,-} \subset (\sigma_n P_i)^{x_1,-}.\]
In particular,
\begin{enumerate}
    \item [(c)] there exists $g \in P_{i-1}$ such that $\sigma_n P_i$ separates $g$ and $x_1C$.
\end{enumerate}
Besides, we have the following facts:
\begin{enumerate}
    \item [(d)] $x_1C$ and $x_nC$ lie on the same side of $\sigma_n P_i$ (by Lemma~\ref{lem-Bsigma}(2));
    \item [(e)] $P_{i-1}$ separates $x_1C$ and $x_nC$ (by Lemma~\ref{lem-Hn-notin-Qn-1}(3)).
\end{enumerate}
By Lemma~\ref{lem-isolateintersect}, we deduce from the above facts (c), (d) and (e) that $P_{i-1} \dc \sigma_n P_i \ne \varnothing$ which contradicts the choice of $P_{i-1}$.
Therefore, $P_i \not\in \mathfrak{B}_\sigma$.
Since $P_i$ is arbitrary in $\mathfrak{A}_i$, it holds that $\mathfrak{A}_i \cap \mathfrak{B}_\sigma = \varnothing$ and we are done.
\end{proof}

\begin{figure}[ht]
    \centering
    \begin{tikzpicture}
        \draw (-3,0) -- (6,0);
        \draw [dashed] (120:3) -- (300:3);
        \path [draw, name path = Piline] (60:3) -- (240:3);
        \path [draw, name path = Pmline] (-2.5,-2.2) -- (5.5,1.5);
        \fill [name intersections = {of = Piline and Pmline, by = {[label=90:$g$]g}}] (g) circle (2pt);
        \coordinate [label=right:$P_{i-1}$] (Pm) at (5.5,1.5);
        \coordinate [label=right:$H_n$] (Hn) at (6,0);
        \coordinate [label=right:$P_i \in \mathfrak{B}_\sigma$] (P) at (60:3);
        \coordinate [label=right:$\sigma_nP_i$] (sP) at (300:3);

        \fill [pattern=dots, pattern color=lightgray] ($(270:2) + (90:0.6)$) -- ($(270:2) + (210:0.6)$) -- ($(270:2) + (330:0.6)$) -- cycle;
        \fill [pattern=dots, pattern color=lightgray] ($(5,0.3) + (90:0.6)$) -- ($(5,0.3) + (210:0.6)$) -- ($(5,0.3) + (330:0.6)$) -- cycle;
        \fill [pattern=dots, pattern color=lightgray] ($(90:2) + (270:0.6)$) -- ($(90:2) + (150:0.6)$) -- ($(90:2) + (30:0.6)$) -- cycle;
        \fill [pattern=dots, pattern color=lightgray] ($(210:2) + (270:0.6)$) -- ($(210:2) + (150:0.6)$) -- ($(210:2) + (30:0.6)$) -- cycle;

        \node (en) at (90:2) {$x_1C$};
        \node (sx1) at (210:2) {$e_{n-1}C$};
        \node (en-1) at (270:2) {$\sigma_nx_1C$};
        \node (xn) at (5,0.3) {$x_nC$};
    \end{tikzpicture}
    \caption{Illustration for the proof of Proposition~\ref{prop-Bsigma-in-B}, supposing $P_i \in \mathfrak{B}_\sigma$.}
    \label{fig-Pi-Bsigma}
\end{figure}

\section{Proof of Theorem~\ref{thm-intersectbound}} \label{sec-bound-intsect}

We employ Ramsey's theorem to simplify the problem.
Let $\Gamma$ be a complete graph on the vertex set $\mathcal{V}$.
By definition, this means each pair of distinct vertices is connected by a unique undirected edge.
An $m$-coloring of $\Gamma$ is a map from the set $\mathcal{E}$ of edges to a set
$\{c_1, c_2, \dots, c_m\}$ of $m$ elements.
The following result is (a special case of) the famous Ramsey's theorem; see, for example, \cite[Section 1.1]{GRS-ramsey}.

\begin{Thm} [Ramsey's theorem] \label{thm-Ramsey}
  For any positive integers $l$ and $m$, there exists a number $R(l;m)$ such that for any complete graph $\Gamma$ with $|\mathcal{V}| \ge R(l;m)$ and any $m$-coloring of $\Gamma$, there exist $i \in \{1, \dots, m\}$ and a subset $\mathcal{V}' \subset \mathcal{V}$ with $|\mathcal{V}'| = l$ such that each edge in $\mathcal{V}'$ (that is, both of the endpoints lie in $\mathcal{V}'$) is colored $c_i$.
\end{Thm}

Let $\mathfrak{Q}$ be an intersecting set.
Let $\Gamma$ be the complete graph on the set $\mathfrak{Q}$, that is, each hyperplane in $\mathfrak{Q}$ is regarded as a vertex of $\Gamma$.
For each distinct pair $P,Q \in \mathfrak{Q}$, we have $ |\operatorname{B}(\alpha_P, \alpha_Q) |< 1$ by Lemma~\ref{lem-finite-dihedral}.
Recall that the set
\[\mathcal{C}ol := \{c \in \mathbb{R} \mid |c|<1, c = \operatorname{B}(\alpha, \beta) \text{ for some } \alpha, \beta \in \Phi^+\}\]
is finite (see Lemma~\ref{lem-cos-finite}), say, $\mathcal{C}ol = \{c_1, \dots, c_m\}$.
In particular, $\operatorname{B}(\alpha_P, \alpha_Q) = c_i$ for some $i$.
Let the edge connecting $P$ and $Q$ be colored $c_i$.
This gives an $m$-coloring of $\Gamma$.

We claim that $|\mathfrak{Q}| < R(\dim V + 2;m)$ (recall that $V$ is the geometric representation space defined in Section~\ref{sec-preliminaries}).
Suppose otherwise $|\mathfrak{Q}| \ge R(\dim V + 2;m)$.
Then, by Ramsey's Theorem~\ref{thm-Ramsey}, there exist $a \in \mathcal{C}ol$ and a subset $\mathfrak{Q}' \subset \mathfrak{Q}$ with $|\mathfrak{Q}'| = \dim V + 2$ such that $\operatorname{B}(\alpha_P, \alpha_Q) = a$ for any $P,Q \in \mathfrak{Q}'$.
Consider the Gram matrix $G := \bigl( \operatorname{B}(\alpha_P, \alpha_Q) \bigr)_{P,Q \in \mathfrak{Q}'}$ which is of the form
\[G = \begin{pmatrix}
    1 & a & \dots & a \\
    a & 1 &  & \vdots \\
    \vdots & & \ddots & a \\
    a & \dots & a & 1
\end{pmatrix}.\]
It is a simple exercise in linear algebra that
\[\operatorname{Rank} G = \begin{cases}
    |\mathfrak{Q}'|-1, & \text{ if } a = -\frac{1}{|\mathfrak{Q}'|-1}, \\
    |\mathfrak{Q}'|, & \text{ otherwise (note that $a \ne 1$)},
\end{cases}\]
and moreover,
\[\operatorname{Rank}\{\alpha_P \mid P \in \mathfrak{Q}'\} \ge \operatorname{Rank} G \ge |\mathfrak{Q}'|-1 = \dim V + 1.\]
But this is absurd since $\{\alpha_P \mid P \in \mathfrak{Q}'\} \subset V$.

In conclusion, we have $|\mathfrak{Q}| < R(\dim V + 2;m)$ for any intersecting set $\mathfrak{Q}$ and hence Theorem~\ref{thm-intersectbound} is proved (note that the numbers $\dim V = |S|$ and $m = |\mathcal{C}ol|$ depend only on the Coxeter group $(W,S)$).

\section{Special cases and examples} \label{sec-eg}

In this section, we assume that $L=\ell$  and denote $N(W)$ for $N_\ell(W)$, that is,
\[N(W) :=\max\{\ell(w_I)\mid I\subset S, \ W_I~{\rm is~finite}\}.\]
Theorem~\ref{thm-exist-intersecting} and Theorem~\ref{thm-intersectbound} imply that the number
\[N'(W) :=\max\{|\mathfrak{Q}| \mid \mathfrak{Q} \text{ is intersecting}\}\]
is finite and that $N'(W)$ is a bound for $(W,S)$. In particular, $N'(W)$ is an upper bound of $\boldsymbol{a}$-function. It is known that any upper bound of $\boldsymbol{a}$-function is $\ge N(W)$. It follows that $N'(W)\ge N(W)$.
We will briefly explain why $N'(W)\le N(W)$ (and hence $N'(W)=N(W)$) in some cases (finite/affine/rank 3) and some details are omitted. However, there also exists an example where $N'(W) > N(W)$.

Let $\mathfrak{Q}$ be an intersecting set of $W$.

\subsection{Finite type}\label{subsec-fin}
If $(W,S)$ is of finite type, then $|\mathfrak{P}|= |\Phi^+| = N(W)$, and hence $|\mathfrak{Q}|\le N(W)$. It follows that $N'(W)\le N(W)$.

\subsection{Affine type}\label{subsec-aff}

Suppose now $(W,S)$ is an affine Weyl group.
Let $\Phi_0$ be the corresponding crystallographic, reduced, finite root system in the sense of \cite{Bourbaki-Lie456}, and $E$ be the real vector space where the dual root system $\Phi_0\spcheck$ lives.
Following \cite[Chapter VI, Section 2]{Bourbaki-Lie456}, for $\alpha \in \Phi_0^+\subset E^*$ and $k \in \mathbb{Z}$, we define the affine hyperplane $L_{\alpha, k}$ of $E$ by
\[L_{\alpha, k} := \{x \in E \mid \langle\alpha,x\rangle = k\}.\]
Let
\[\mathscr{H} = \{L_{\alpha, k} \mid \alpha \in \Phi^+_0, k \in \mathbb{Z}\}.\]
It is known that $\mathscr{H}$ is in one-to-one correspondence with the set $\mathfrak{P}$ of hyperplanes in $V^*$; see, for example, \cite[Section~6.5]{Hum}.
Moreover, for two distinct hyperplanes $P_1, P_2 \in \mathfrak{P}$, $P_1 \dc P_2 \ne \varnothing$ if and only if their corresponding hyperplanes in $\mathscr{H}$ are of different directions (we say that $L_{\alpha, k}$ is of direction $\alpha$).
It follows that $|\mathfrak{Q}|\le|\Phi_0^+|=N(W)$ and hence $N'(W)\le N(W)$.

\subsection{Hyperbolic type of rank 3}\label{subsec-r3}

If $(W,S)$ is of rank 3 and $(W,S)$ is not of finite types $A_3, B_3,H_3$ or affine types $\widetilde{A}_2, \widetilde{B}_2, \widetilde{G}_2$, then $(W,S)$ is of hyperbolic type (see \cite[Section 6.7]{Hum}).

Let $(W,S)$ be a hyperbolic Coxeter group with $S=\{s_1,s_2,s_3\}$, and $m, n, p$ be the order of $s_1s_2, \, s_2s_3, \, s_1s_3$, respectively ($m,n,p\in\{2,3,\dots\}\cup\{\infty\}$). Then the bilinear form $\operatorname{B}(-,-)$ on $V$ has signature $(2,1)$. One identifies $V^*$ with $V$ by means of $B(-,-)$. Let $\nu$ be the canonical projection from $V\backslash\{0\}$ to the projective space $\mathbb{P}V$. Let $N=\{v\in V\mid \operatorname{B}(v,v)<0\}$ and $\mathbb{H}^2=\nu(N)$ (called hyperbolic plane). By a line in $\mathbb{H}^2$ we mean a set of the form $\nu(P\cap N)$, where $P$ is a hyperplane in $V$ such that $P\cap N\ne \varnothing$. One can define angles in $\mathbb{H}^2$; see \cite[5.4]{And05}.

We have the following facts.
\begin{enumerate}
    \item The number $N(W)$ equals the maximal number in $\{m,n,p\}\backslash\{\infty\}$.
    \item Each hyperplane in $\mathfrak{P}$ is identified with a line in $\mathbb{H}^2$ (Let $P\in\mathfrak{P}$. Since $U^\circ\subset N$ by \cite[Exercise \S 4 (12)(a)]{Bourbaki-Lie456} and $P\cap U^\circ\ne \varnothing$ by Corollary \ref{cor-single-intersecting}, it follows that $P\cap N\ne \varnothing$. The plane $P$ is identified with the image of $P\cap N$ in $\mathbb{P}V$, which is a line in $\mathbb{H}^2$).
    \item If $P_1,P_2\in\mathfrak{P}$, $P_1\ne P_2$, and $P_1\dc P_2\ne \varnothing$, then the lines corresponding to $P_1,P_2$ intersect at a point in $\mathbb{H}^2$ (because $U^\circ\subset N$).
    \item The angles between two different intersecting lines (arising from $\mathfrak{P}$) are in $A=\{\frac{k\pi}{N}\mid N\in\{m,n,p\}\backslash\{\infty\},0<k<N\}$ (because any finite subgroup of $W$ is conjugate to a subgroup of a finite standard parabolic subgroup of $W$; see \cite[Theorem~4.5.3]{BB05}).
    \item The angle sum of a triangle in  $\mathbb{H}^2$ is less than $\pi$ (see \cite[5.4]{And05}).
\end{enumerate}
Due to (2), we use same notation for a hyperplane in $\mathfrak{P}$ and its corresponding line in $\mathbb{H}^2$.

Let $\mathfrak{Q}=\{P_0,P_1,\cdots,P_{t-1}\}$, viewed as a set of $t$ lines that intersect pairwise in $\mathbb{H}^2$,
and denote $a_i=P_i\cap P_0$, $(1\le i\le t-1)$. Let $\Lambda,\Lambda'$ be two connected components of $\mathbb{H}^2\backslash P_0$ and choose $b_i\in P_i\cap\Lambda$ $(1\le i\le t-1)$. There is an $a\in P_0$ such that all $a_i$ lie on the same side of $a$.

With the above preparations, we define the hyperbolic angles $\theta_i=\angle b_ia_ia$ $(1\le i\le t-1)$. We claim that
\begin{equation}\label{mingap}
    |\theta_i-\theta_j|\ge\frac{\pi}{N(W)},~i\ne j.
\end{equation}
To see this, let $a_{ij} :=P_i\cap P_j$,  and we consider the following cases on the relative position of lines $P_0,P_i,P_j$ in $\mathbb{H}^2$.
\begin{enumerate}
    \item [(a)] $P_0\cap P_i\cap P_j\ne \varnothing$, that is, $a_{ij} = a_i = a_j$. In this case, it is clear that $|\theta_i-\theta_j| = \angle b_i a_{ij} b_j \ge \frac{\pi}{N(W)}$ by the facts (1) and (4) (see Figure~\ref{fig-a}).
\begin{figure}[ht]
    \centering
    \begin{tikzpicture}
        \draw (-2,0) -- (3,0);
        \draw (60:3) -- (240:1);
        \draw (120:3) -- (300:1);
        \coordinate [label=left:$P_0$] (Hn) at (-2,0);
        \coordinate [label=right:$P_j$] (P) at (60:3);
        \coordinate [label=left:$P_i$] (sP) at (120:3);
        \coordinate [label=above left:$a_{ij}$] (sP) at (-0.2,0);
        \coordinate [label=above:$\theta_i$] (sP) at (0,0.1);
        \coordinate [label=above right:$\theta_j$] (sP) at (0.4,0);
        \filldraw[fill=black] (2.5,0) circle[radius=0.3mm] node [above] {$a$};
        \filldraw[fill=black] (60:2) circle[radius=0.3mm] node [below right] {$b_j$};
        \filldraw[fill=black] (120:2) circle[radius=0.3mm] node [below left] {$b_i$};
        \draw (0,0) -- +(0:0.5) arc (0:60:0.5);
        \draw (0,0) -- +(0:0.2) arc (0:120:0.2);
        \node (or) at (3.3,1) {or};
    \end{tikzpicture}
    \begin{tikzpicture}
        \draw (-2,0) -- (3,0);
        \draw (60:3) -- (240:1);
        \draw (120:3) -- (300:1);
        \coordinate [label=left:$P_0$] (Hn) at (-2,0);
        \coordinate [label=right:$P_i$] (P) at (60:3);
        \coordinate [label=left:$P_j$] (sP) at (120:3);
        \coordinate [label=above left:$a_{ij}$] (sP) at (-0.2,0);
        \coordinate [label=above:$\theta_j$] (sP) at (0,0.1);
        \coordinate [label=above right:$\theta_i$] (sP) at (0.4,0);
        \filldraw[fill=black] (2.5,0) circle[radius=0.3mm] node [above] {$a$};
        \filldraw[fill=black] (60:2) circle[radius=0.3mm] node [below right] {$b_i$};
        \filldraw[fill=black] (120:2) circle[radius=0.3mm] node [below left] {$b_j$};
        \draw (0,0) -- +(0:0.5) arc (0:60:0.5);
        \draw (0,0) -- +(0:0.2) arc (0:120:0.2);
    \end{tikzpicture}
    \caption{Illustration for case (a)}
    \label{fig-a}
\end{figure}

\item [(b)] $P_i\cap P_j\in\Lambda$. Applying (5) to the triangle $a_{ij}a_ia_j$ yields $\theta_i+ \angle a_i a_{ij} a_j <\theta_j$ or $\theta_j+ \angle a_i a_{ij} a_j <\theta_i$ (see Figure~\ref{fig-b}). Moreover, the facts (1) and (4) imply that $\angle a_i a_{ij} a_j \ge\frac{\pi}{N(W)}$. Thus, \eqref{mingap} follows immediately.

\begin{figure}[ht]
\begin{tikzpicture}
\path [draw, name path=Lline] (-2,0) -- (3.5,0);
\path [draw, name path=Ljline] (-0.5,2.59808) -- (1.5,-0.866025);
\path [draw, name path=Liline] (0.5,2.59808) -- (-1.5,-0.866025);
\coordinate [label=left:$P_j$] (lj) at (-0.5,2.59808);
\coordinate [label=right:$P_i$] (li) at (0.5,2.59808);
\coordinate [label=left:$P_0$] (l) at (-2,0);
\fill [name intersections = {of = Liline and Ljline, by = {[label=left:$a_{ij}$]aij}}] (aij) circle (1pt);
\fill [name intersections = {of = Lline and Liline, by = {[label=above left:$a_i$]ai}}] (ai) circle (1pt);
\fill [name intersections = {of = Lline and Ljline, by = {[label=below left:$a_j$]aj}}] (aj) circle (1pt);
\draw (-1,0) -- +(0:0.15) arc (0:60:0.15);
\coordinate [label=above right:$\theta_i$] (thetai) at (-0.9,-0.1);
\draw (1,0) -- +(0:0.15) arc (0:120:0.15);
\coordinate [label=above right:$\theta_j$] (thetaj) at (1,-0.1);
\coordinate [label=above right:$\Lambda$] (Lambda) at (3,1);
\coordinate [label=above right:$\Lambda'$] (Lambdap) at (3,-1);
\filldraw[fill=black] (3,0) circle[radius=0.3mm] node [above] {$a$};
\filldraw[fill=black] (-0.43,1) circle[radius=0.3mm] node [left] {$b_i$};
\filldraw[fill=black] (0.43,1) circle[radius=0.3mm] node [right] {$b_j$};
\node (or) at (4,1) {or};
\end{tikzpicture}
\begin{tikzpicture}
\path [draw, name path=Lline] (-2,0) -- (3.5,0);
\path [draw, name path=Ljline] (-0.5,2.59808) -- (1.5,-0.866025);
\path [draw, name path=Liline] (0.5,2.59808) -- (-1.5,-0.866025);
\coordinate [label=left:$P_i$] (lj) at (-0.5,2.59808);
\coordinate [label=right:$P_j$] (li) at (0.5,2.59808);
\coordinate [label=left:$P_0$] (l) at (-2,0);
\fill [name intersections = {of = Liline and Ljline, by = {[label=left:$a_{ij}$]aij}}] (aij) circle (1pt);
\fill [name intersections = {of = Lline and Liline, by = {[label=above left:$a_j$]ai}}] (ai) circle (1pt);
\fill [name intersections = {of = Lline and Ljline, by = {[label=below left:$a_i$]aj}}] (aj) circle (1pt);
\draw (-1,0) -- +(0:0.15) arc (0:60:0.15);
\coordinate [label=above right:$\theta_j$] (thetai) at (-0.9,-0.1);
\draw (1,0) -- +(0:0.15) arc (0:120:0.15);
\coordinate [label=above right:$\theta_i$] (thetaj) at (1,-0.1);
\coordinate [label=above right:$\Lambda$] (Lambda) at (3,1);
\coordinate [label=above right:$\Lambda'$] (Lambdap) at (3,-1);
\filldraw[fill=black] (3,0) circle[radius=0.3mm] node [above] {$a$};
\filldraw[fill=black] (-0.43,1) circle[radius=0.3mm] node [left] {$b_j$};
\filldraw[fill=black] (0.43,1) circle[radius=0.3mm] node [right] {$b_i$};
\end{tikzpicture}
\caption{Illustration for case (b)} \label{fig-b}
\end{figure}

\item [(c)] $P_i\cap P_j\in\Lambda'$. Using similar arguments as in (b) (see Figure~\ref{fig-c}), we conclude that \eqref{mingap} holds.

\begin{figure}[ht]
\begin{tikzpicture}
\path [draw, name path=Lline] (-2,0) -- (3.5,0);
\path [draw, name path=Liline] (-1.5,0.866025) -- (0.5,-2.59808);
\path [draw, name path=Ljline] (1.5,0.866025) -- (-0.5,-2.59808);
\coordinate [label=left:$P_0$] (l) at (-2,0);
\coordinate [label=left:$P_i$] (li) at (-1.5,0.866025);
\coordinate [label=right:$P_j$] (lj) at (1.5,0.866025);
\fill [name intersections = {of = Liline and Ljline, by = {[label=left:$a_{ij}$]aij}}] (aij) circle (1pt);
\fill [name intersections = {of = Lline and Liline, by = {[label=below left:$a_i$]ai}}] (ai) circle (1pt);
\fill [name intersections = {of = Lline and Ljline, by = {[label=below right:$a_j$]aj}}] (aj) circle (1pt);
\coordinate [label=right:$\Lambda$] (Lambda) at (3,1);
\coordinate [label=right:$\Lambda'$] (Lambdap) at (3,-1);
\draw (-1,0) -- +(0:0.15) arc (0:120:0.15);
\draw (1,0) -- +(0:0.15) arc (0:60:0.15);
\draw (1,0) -- +(180:0.15) arc (180:240:0.15);
\draw (1,0) -- +(0:0.2) arc (0:60:0.2);
\draw (1,0) -- +(180:0.2) arc (180:240:0.2);
\coordinate [label=above right:$\theta_i$] (thetai) at (-1,0);
\coordinate [label=above right:$\theta_j$] (thetaj) at (1.2,-0.1);
\coordinate [label=below left:$\theta_j$] (thetaj) at (0.9,0.1);
\filldraw[fill=black] (3,0) circle[radius=0.3mm] node [above] {$a$};
\filldraw[fill=black] (-1.3,0.5) circle[radius=0.3mm] node [above right] {$b_i$};
\filldraw[fill=black] (1.3,0.5) circle[radius=0.3mm] node [above left] {$b_j$};
\node (or) at (4,-1) {or};
\end{tikzpicture}
\begin{tikzpicture}
\path [draw, name path=Lline] (-2,0) -- (3.5,0);
\path [draw, name path=Liline] (-1.5,0.866025) -- (0.5,-2.59808);
\path [draw, name path=Ljline] (1.5,0.866025) -- (-0.5,-2.59808);
\coordinate [label=left:$P_0$] (l) at (-2,0);
\coordinate [label=left:$P_j$] (li) at (-1.5,0.866025);
\coordinate [label=right:$P_i$] (lj) at (1.5,0.866025);
\fill [name intersections = {of = Liline and Ljline, by = {[label=left:$a_{ij}$]aij}}] (aij) circle (1pt);
\fill [name intersections = {of = Lline and Liline, by = {[label=below left:$a_j$]ai}}] (ai) circle (1pt);
\fill [name intersections = {of = Lline and Ljline, by = {[label=below right:$a_i$]aj}}] (aj) circle (1pt);
\coordinate [label=right:$\Lambda$] (Lambda) at (3,1);
\coordinate [label=right:$\Lambda'$] (Lambdap) at (3,-1);
\draw (-1,0) -- +(0:0.15) arc (0:120:0.15);
\draw (1,0) -- +(0:0.15) arc (0:60:0.15);
\draw (1,0) -- +(180:0.15) arc (180:240:0.15);
\draw (1,0) -- +(0:0.2) arc (0:60:0.2);
\draw (1,0) -- +(180:0.2) arc (180:240:0.2);
\coordinate [label=above right:$\theta_j$] (thetai) at (-1,-0.1);
\coordinate [label=above right:$\theta_i$] (thetaj) at (1.2,-0.1);
\coordinate [label=below left:$\theta_i$] (thetaj) at (0.9,0.1);
\filldraw[fill=black] (3,0) circle[radius=0.3mm] node [above] {$a$};
\filldraw[fill=black] (-1.3,0.5) circle[radius=0.3mm] node [above right] {$b_j$};
\filldraw[fill=black] (1.3,0.5) circle[radius=0.3mm] node [above left] {$b_i$};
\end{tikzpicture}
\caption{Illustration for case (c)} \label{fig-c}
\end{figure}
\end{enumerate}

Due to \eqref{mingap}, up to a rearrangement, one can assume that $0<\theta_1<\cdots<\theta_{t-1}<\pi$ and $\theta_{i+1}-\theta_i\ge\frac{\pi}{N(W)}$ for $1\le i<t-1$. In particular, we have $\frac{(t-1)\pi}{N(W)}\le\theta_{t-1}<\pi$. It follows that $|\mathfrak{Q}|=t\le N(W)$, and hence $N'(W)\le N(W)$.

\begin{Rem}
In the case of hyperbolic type of rank 3, since the equality in \eqref{mingap} may hold only in Case (a), it follows that if $|\mathfrak{Q}|=N(W)$, then the lines in $\mathfrak{Q}$ must be concurrent. Clearly, this is not true even for $\widetilde{A}_2$ (consider, for example, 3 lines forming the fundamental alcove).
\end{Rem}

\begin{Rem}
  After more detailed verifications, our approach is capable to give the conjectural bound $N_L(W)$ in the cases finite/affine/rank 3 for arbitrary weight function $L$.
  Since this bound have been verified in these cases, we omit the details here to avoid making this section too lengthy.
\end{Rem}

\subsection{An example for \texorpdfstring{$N'(W) > N(W)$}{N'>N}}

As shown in the following example, for a certain Coxeter group $(W,S)$ there exists an intersecting set $\mathfrak{Q}$ such that $|\mathfrak{Q}| > N(W)$.
Thus, it holds $N'(W)>N(W)$ for this group $W$.

\begin{Ex} \label{ex-greater-bound}
Consider the Coxeter group $(W,S)$ with Coxeter graph as follows:
$$
\begin{tikzpicture}
\filldraw[fill=black] (0,-0.75) circle[radius=0.5mm] node [above] {$s_2$};
\filldraw[fill=black] (-1.5,-0.75) circle[radius=0.5mm] node [left] {$s_1$};
\filldraw[fill=black] (1.29904,0) circle[radius=0.5mm] node [right] {$s_3$};
\filldraw[fill=black] (1.29904,-1.5) circle[radius=0.5mm] node [right] {$s_4$};
\draw (-1.5,-0.75) -- (0,-0.75) node [above left] {$6~\quad$};
\draw (0,-0.75) -- (1.29904,0);
\draw (0,-0.75) -- (1.29904,-1.5);
\end{tikzpicture}
$$
Let $I=\{s_2,s_3,s_4\}$. Then $N(W)=\ell(w_I)=6$. All positive roots of $W_I$ are $\beta_1=\alpha_3$, $\beta_2=\alpha_2$, $\beta_3=\alpha_4$, $\beta_4=\alpha_3+\alpha_2$, $\beta_5=\alpha_2+\alpha_4$, $\beta_6=\alpha_3+\alpha_2+\alpha_4$. Since $W_I$ is finite, $\{H_{\beta_i}\mid1\le i\le 6\}$ is intersecting. Let $b_i=\operatorname{B}(\alpha_1,\beta_i)$ $(1\le i\le6)$. Simple calculations show that $b_1=b_3=0$ and $b_2=b_4=b_5=b_6=-\frac{\sqrt{3}}{2}$. It follows that $H_{\alpha_1}\dc H_{\beta_i}\ne \varnothing$ $(1\le i\le 6)$ by Lemma~\ref{lem-finite-dihedral}. Thus, $\{H_{\alpha_1},H_{\beta_1},\dots,H_{\beta_6}\}$ is intersecting.
\end{Ex}

\begin{Rem}
    Let $W_I$ be a parabolic subgroup such that $\ell(w_I)=N(W)$, and $\mathfrak{P}_I\subset\mathfrak{P}$ be the set of hyperplanes corresponding to roots of $W_I$.  Example \ref{ex-greater-bound} also indicates that there might be a hyperplane $P$ in $\mathfrak{P}\backslash\mathfrak{P}_I$ that intersects (in $U^\circ$) all hyperplanes in $\mathfrak{P}_I$. However, it is clear that such $P$ does not exist in affine case.
\end{Rem}

\appendix

\section{Paths in Bruhat graph}

This section is a reformulation of \eqref{eq-expansion-TxTy} as well as condition \eqref{eq-<} in terms of paths in the Bruhat graph, due to an anonymous referee.

Let $T := \bigcup_{w \in W} wSw^{-1}$ be the set of reflections.
As in Section~\ref{subsec-overview-proof}, let $x, y \in W$ and we fix a reduced expression $y = s_1 s_2 \cdots s_k$ for $y$.
Let \[t_i := s_1 \cdots s_{i-1} s_i s_{i-1} \cdots s_1 \in T \quad (i = 1, \dots, k). \]
The condition \eqref{eq-<} (for all $n = 1,2, \dots, p$) is equivalent to the following path in Bruhat order,
\begin{equation}\label{eq-reform-TxTy}
P: \ xy < x t_{i_1}y < x t_{i_1}t_{i_2}y < \dots < x t_{i_1} t_{i_2} \cdots t_{i_p} y.
\end{equation}
To see this, we need the following lemma.
\begin{Lem}  \label{lem-app1}
    Suppose $t \in T$ is a reflection, $u, u' \in W$, and $tu > u$.
    Then $u't< u'$ if and only if $u'tu < u'u$.
\end{Lem}

\begin{proof}
     It is a standard fact that for $w \in W$ and $r \in T$, $wr < w$ if and only if $w \alpha_r \in \Phi^-$, where $\alpha_r$ is the positive root corresponding to $r$ (see \cite[Proposition 5.7]{Hum}).
     Rewrite $u'tu$ as $u'u (u^{-1}tu)$.
     Then, $\alpha_{u^{-1}tu} = u^{-1} \alpha_t \in \Phi^+$ and $u'u \alpha_{u^{-1}tu} = u'\alpha_t$.
     Therefore, $u'tu < u'u$ if and only if $u'\alpha_t \in \Phi^-$, which is equivalent to $u't < u'$.
\end{proof}

For each $n = 1, \dots, p$, we set in Lemma \ref{lem-app1}
\[u' := xs_1 s_2\cdots \widehat{s_{i_1}}\cdots\widehat{s_{i_{n-1}}}\cdots s_{i_n-1}, \ t:= s_{i_n}, \text{ and } u :=  s_{i_n+1} s_{i_n+2} \cdots s_k,\]
then we see that $u't < u'$ if and only if $x t_{i_1} \cdots t_{i_{n-1}}y < xt_{i_1} \cdots t_{i_n} y$.
This interprets condition \eqref{eq-<} (for all $n=1,\dots, n$) into the path $P$ in \eqref{eq-reform-TxTy}.

We define $L'(w s w^{-1}) := L(s)$ for all $s \in S$ and $w \in W$.
This is a well-defined function on $T$.
Then \eqref{eq-expansion-TxTy} can be written equivalently as
\[T_x T_y = \sum_P \zeta_P T_{z_P}\]
where $P$ runs over all the paths in Bruhat order of the form \eqref{eq-reform-TxTy}, $z_P := x t_{i_1} t_{i_2} \cdots t_{i_p} y$, and $\zeta_P := (v^{L'(t_{i_1})} - v^{-L'(t_{i_1})} )  (v^{L'(t_{i_2})} - v^{-L'(t_{i_2})} ) \cdots  (v^{L'(t_{i_p})} - v^{-L'(t_{i_p})} )$.

\section{A reformulation of Theorem~\ref{thm-exist-intersecting}}

Based on the comments from another anonymous referee, this section consists of a reformulation of Theorem~\ref{thm-exist-intersecting} and its proof.
The proof, mostly given by the referee, uses other techniques and concepts such as the dominance order in the root system.
The reader might find it more helpful to understand our results.

Let $\xi_s :=v^{L(s)}-v^{-L(s)}$ $(s\in S)$.
Given $x, y\in W$, we write (compared with \eqref{eq-expansion-TxTy}) that
\begin{equation} \label{eq-expansion-TxTy-app}
    T_xT_y = \sum\nolimits_{z \in W} f_z T_z
\end{equation}
where $f_z$ is a polynomial in variables $\{\xi_s \mid s \in S\}$.
We regard each $\xi_s$ to be of degree one and define the degree of $f_z$ to be its degree as a polynomial in those $\xi_s$.

A {\bfseries \itshape 2-spherical set} is defined to be a subset $\mathcal{Q}$ of $\Phi^+$ with $|{\rm B}(\alpha,\beta)|<1$ for any distinct roots $\alpha,\beta\in\mathcal{Q}$.
In view of Lemma \ref{lem-finite-dihedral}, this is equivalent to saying that the corresponding hyperplanes $H_\alpha$ and $H_\beta$ are intersecting.
By convention, the empty set $\varnothing$ and singletons $\{\alpha\}$ $(\alpha\in\Phi^+)$ are regarded to be 2-spherical.

For any $w\in W$, set $N(w)=\{\alpha\in\Phi^+\mid w\alpha\in\Phi^-\}$.
The following is a reformulation of Theorem \ref{thm-exist-intersecting}.
\begin{Thm}\label{exist-2-sphere-set}
    If $f_zT_z$ is a term in \eqref{eq-expansion-TxTy-app}, then there is a 2-spherical set $\mathcal{Q}_z \subset N(y)\cap N(z)$ with $|\mathcal{Q}_z|= \deg f_z$.
\end{Thm}

Let us outline its proof, which is by induction on $\ell(y)$.
If $y=e$, we just take $\mathcal{Q}_z=\varnothing$.
For the case $y = s \in S$ it is also easy to see that $\mathcal{Q}_z = \varnothing$ or $\{\alpha_s\}$ suffices.

Suppose that each term $f_zT_z$ in \eqref{eq-expansion-TxTy-app} gives a
2-spherical set $\mathcal{Q}_z\subset N(y)\cap N(z)$ with $|\mathcal{Q}_z|= \deg f_z$.
Suppose $s\in S$ and $ys>y$ (hence $\alpha_s \notin N(y)$).
We have
\[T_xT_{ys} = \sum\nolimits_{z \in W} f_z T_z T_s = \sum\nolimits_{z\in W} f_z T_{zs} + \sum\nolimits_{z\in W, \, zs < z} f_z \xi_s T_z.\]

\textbf{(I).} For the term $f_z T_{zs}$, we need to find a 2-spherical set $\mathcal{Q}'$ such that $|\mathcal{Q}'|= |\mathcal{Q}_z|$ and $\mathcal{Q}'\subset N(ys)\cap N(zs)$.
It turns out that $\mathcal{Q}'=s\mathcal{Q}_z$ suffices (we leave the details to the reader).

\textbf{(II).} For the term $f_z \xi_s T_{z}$ ($zs < z$), we need to find a set $\mathcal{Q}'\subset \Phi^+$ such that
\begin{center}
    (i). $\mathcal{Q}'$ is 2-spherical; \quad
     (ii). $|\mathcal{Q}'|= |\mathcal{Q}_z| + 1$; \quad and
    (iii). $\mathcal{Q}'\subset N(ys)\cap N(z)$.
\end{center}
The construction of such a $\mathcal{Q}'$ is parallel to that in Theorem~\ref{thm-exist-intersecting-2}.
We will partition $\mathcal{Q}_z$ as $\mathcal{Q}_z = \mathcal{A} \sqcup \mathcal{B}$ and define
\begin{equation}\label{eq-Q'}
    \mathcal{Q}' :=\mathcal{A}\cup s\mathcal{B}\cup\{\alpha_s\}.
\end{equation}
Let
$\mathcal{A}_0:=\{\alpha\in\mathcal{Q}_z\mid s\alpha \notin N(z)\}$.
Since we require $\mathcal{Q}' \subset N(z)$, $\mathcal{A}_0$ must be contained in $\mathcal{A}$.
For $i \ge 1$ we define inductively
\[\mathcal{A}_{i}:=\{\gamma\in\mathcal{Q}_z\setminus(\mathcal{A}_0\sqcup\mathcal{A}_1\sqcup\cdots\sqcup\mathcal{A}_{i-1}) \mid |\operatorname{B}(s\gamma,\beta)|\ge1~{\rm for~some~}\beta\in\mathcal{A}_{i-1}\}.\]
Since the set $\mathcal{Q}'$ is required to be 2-spherical, all the sets $\mathcal{A}_i$ must be contained in $\mathcal{A}$ as well.
Clearly, there is a smallest $l$ such that $\mathcal{A}_i=\varnothing$ for any $i > l$.
Set $\mathcal{A}:=\bigsqcup_{i=0}^{l}\mathcal{A}_i$ and $\mathcal{B}:=\mathcal{Q}_z\backslash\mathcal{A}$.
Finally, define $\mathcal{Q}'$ as in \eqref{eq-Q'}.
We need to verify the properties (i), (ii) and (iii) above.

For $\alpha, \beta \in \Phi$, we say that $\alpha$ {\bfseries \textit{dominates}} $\beta$ (and write $\alpha \operatorname{dom} \beta$) if $w\beta \in \Phi^-$ whenever $w \alpha \in \Phi^-$ ($w \in W$).
This gives a partial order on $\Phi$.
See \cite{BH93, Fu12} for more details.

\begin{Lem}\label{Bge1} \leavevmode
    \begin{enumerate}
        \item {\upshape(\cite[Lemma 3.2(ii)]{Fu12})} The dominance order is $W$-invariant, that is, if $\alpha \operatorname{dom} \beta$ then $w \alpha \operatorname{dom} w \beta$ for any $w \in W$.
        \item {\upshape(\cite[Lemma 3.2(v)]{Fu12})} $\operatorname{B}(\alpha,\beta)\ge1$ if and only if $\alpha \operatorname{dom} \beta$ or $\beta \operatorname{dom} \alpha$.
        \item {\upshape(\cite[Lemma 4.7.4]{BB05})} Let $\beta \in \Phi^+$ and $s \in S$.
            Then $\beta \operatorname{dom} \alpha_s$ if and only if ${\rm B}(\beta,\alpha_s)\ge1$.
    \end{enumerate}
\end{Lem}

\begin{Lem}\label{Bge-1}
If $w\in W$ and $\alpha,\beta\in N(w)$, then $\operatorname{B}(\alpha,\beta)>-1$.
\end{Lem}
\begin{proof}
Otherwise, there exists infinitely many positive roots of the form $\lambda\alpha+\mu\beta$ with $\lambda,\mu\ge0$ (see, for example, \cite[Proposition 4.5.4]{BB05}).
These roots are contained in $N(w)$ contradicting the finiteness of $N(w)$.
\end{proof}

\begin{proof}[Proof of {\upshape(i)}]
    The emptiness of $\mathcal{A}_l$ implies that $\mathcal{A}\cup s\mathcal{B}$ is 2-spherical.
    For any $\alpha\in \mathcal{Q}_z$, since $\alpha, \alpha_s \in N(z)$, we have $\operatorname{B}(\alpha,\alpha_s)>-1$ by Lemma~\ref{Bge-1}.
    On the other hand, since $y \alpha \in \Phi^-$ and $y \alpha_s \in \Phi^+$, we deduce $\operatorname{B}(\alpha,\alpha_s) < 1$ by Lemma~\ref{Bge1}(3).
    To conclude, $|\operatorname{B}(\alpha,\alpha_s)|<1$ for any $\alpha\in\mathcal{Q}_z$ (hence for any $\alpha\in s\mathcal{Q}_z$).
    In particular,  $|\operatorname{B}(\alpha,\alpha_s)|<1$ for any $\alpha\in\mathcal{A}\cup s\mathcal{B}$.
    Therefore, $\mathcal{Q}'$ is 2-spherical.
\end{proof}

\begin{proof}[Proof of {\upshape(ii)}]
    It is clear that $\alpha_s \notin \mathcal{A} \cup s\mathcal{B}$ since $\mathcal{A} \subset N(y)$ and $\mathcal{B} \subset \Phi^+$.
    We also have $\mathcal{A}_0 \cap s\mathcal{B} = \varnothing$ by the definition of $\mathcal{A}_0$ and the fact that $\mathcal{B} \subset N(z)$.
    For $i > 0$, suppose $\gamma \in \mathcal{A}_i \cap s\mathcal{B}$.
    By the definition of $\mathcal{A}_i$, there is a $\beta \in \mathcal{A}_{i-1}$ such that $|\operatorname{B}(\beta, s\gamma)| \ge 1$.
    But since $s\gamma \in \mathcal{B}$ this contradicts the fact that $\mathcal{Q}_z$ is 2-spherical.
    In conclusion, $\mathcal{A}\cup s\mathcal{B}\cup\{\alpha_s\}$ is a disjoint union and $|\mathcal{Q}'| = |\mathcal{Q}_z| +1$.
\end{proof}

We divide the proof of (iii) into the following three steps.

\noindent\textit{Step 1.} Show that $\mathcal{A}\subset N(z)$ and $s\mathcal{B}\cup\{\alpha_s\}\subset N(ys)\cap N(z)$.

The first assertion is clear since $\mathcal{A}\subset\mathcal{Q}_z\subset N(z)$.
The fact $\alpha_s\in N(ys)\cap N(z)$ follows easily from Lemma~\ref{lem-chamber}(3).
For any $\beta\in \mathcal{B}$, we have $\beta\in N(y)$, and hence $s\beta\in N(ys)$.
We also have $s\beta \in N(z)$ for otherwise $\beta$ would fall in $\mathcal{A}_0$, contradicting $\beta \in \mathcal{B}$.

\noindent\textit{Step 2.} Show that $\mathcal{A}_0\cup\mathcal{A}_1\subset N(ys)$.

Let $\alpha\in\mathcal{A}_0$.
Since $zs\alpha=z\alpha-2\operatorname{B}(\alpha,\alpha_s)z\alpha_s$ and $z\alpha,z\alpha_s\in\Phi^-$, the condition $zs\alpha\in\Phi^+$ forces
\begin{equation} \label{eq-app-baasge0}
\operatorname{B}(\alpha,\alpha_s)>0 \quad \text{(this holds for all $\alpha \in \mathcal{A}_0$)}.
\end{equation}
This implies that $ys\alpha=y\alpha-2\operatorname{B}(\alpha,\alpha_s)y\alpha_s\in\Phi^-$ since  $y\alpha\in\Phi^-$ and $y\alpha_s\in\Phi^+$.
Therefore, $\mathcal{A}_0\subset N(ys)$.

For any $\beta\in\mathcal{A}_1$, by the definition of $\mathcal{A}_1$, there exists $\alpha\in\mathcal{A}_0$ such that $|\operatorname{B}(\beta,s\alpha)| \ge1$.
In particular, $s\alpha\in N(y)$ since $\mathcal{A}_0\subset N(ys)$.
Since $\beta,s\alpha\in N(y)$, Lemma~\ref{Bge-1} implies $\operatorname{B}(\beta,s\alpha)>-1$, and hence $\operatorname{B}(\beta,s\alpha)\ge1$.
This in turn implies that $\beta \operatorname{dom} s\alpha$ or $s\alpha\operatorname{dom} \beta$ by Lemma~\ref{Bge1}(2).
Since $z\beta\in\Phi^-$ but $zs\alpha\in\Phi^+$, it is not true that $\beta \operatorname{dom} s\alpha$.
Therefore, $s\alpha \operatorname{dom} \beta$, and hence $\alpha \operatorname{dom} s\beta$ by  Lemma~\ref{Bge1}(1).
Since $y \alpha \in \Phi^-$, this implies that $ys\beta\in\Phi^-$, and hence $\mathcal{A}_1\subset N(ys)$.

As a byproduct, since $\mathcal{Q}_z$ is 2-spherical, we have
\begin{equation} \label{eq-app-1le<1-2}
1 \le \operatorname{B}(\beta,s\alpha)= \operatorname{B}(\beta,\alpha) - 2 \operatorname{B}(\alpha,\alpha_s) \operatorname{B}(\beta,\alpha_s) < 1 - 2 \operatorname{B} (\alpha,\alpha_s) \operatorname{B}(\beta,\alpha_s).
\end{equation}
This, together with \eqref{eq-app-baasge0}, implies that
\begin{equation} \label{eq-app-bbasle0}
    \operatorname{B}(\beta,\alpha_s)<0 \quad \text{(this holds for all $\beta \in \mathcal{A}_1$)}.
\end{equation}

\noindent\textit{Step 3.} Show that $\mathcal{A}_i \subset N(ys)$ for all $i$.

For $i\ge1$, $\delta\in\mathcal{A}_i$, suppose $\gamma\in\mathcal{A}_{i-1}$ is chosen such that $|\operatorname{B}(\delta,s\gamma)|\ge1$.
We will prove by induction that
\[\text{$\delta \operatorname{dom} s\gamma$ and $\operatorname{B}(\delta, \alpha_s) > 0$ if $i$  is  even;} \quad \text{$s\gamma \operatorname{dom} \delta$  and $\operatorname{B}(\delta,\alpha_s)<0$  if $i$ is odd,}\]
and $\delta\in  N(ys)$ for any $i$.
The case $i=1$ has been settled in the previous step.

Let $i\ge2$ and $\beta\in\mathcal{A}_{i-2}$ be chosen such that $|\operatorname{B}(s\gamma,\beta)|\ge1$.
By induction, we have $\gamma\in N(ys)$, hence $s\gamma\in N(y)$.
Thus, Lemma \ref{Bge-1} implies that $\operatorname{B}(\delta,s\gamma)>-1$, and hence $\operatorname{B}(\delta,s\gamma)\ge1$ (since $|\operatorname{B}(\delta,s\gamma)|\ge1$).
Thus, Lemma \ref{Bge1}(2) implies that $\delta \operatorname{dom} s\gamma$ or $s\gamma \operatorname{dom} \delta$.

As in \eqref{eq-app-1le<1-2}, we have
$1\le\operatorname{B}(\delta,s\gamma)<1-2\operatorname{B}(\delta,\alpha_s)\operatorname{B}(\gamma,\alpha_s)$,
which implies that $\operatorname{B}(\delta,\alpha_s)\operatorname{B}(\gamma,\alpha_s)<0$.

Suppose that $i$ is even.
By induction, we have $s\beta \operatorname{dom} \gamma$ (equivalently, $\beta \operatorname{dom} s\gamma$) and $\operatorname{B}(\gamma,\alpha_s)<0$, and hence $\operatorname{B}(\delta,\alpha_s)>0$.
Since  $y\delta\in\Phi^-$ and $y\alpha_s\in\Phi^+$, this gives $ys\delta=y\delta-2\operatorname{B}(\delta,\alpha_s)y\alpha_s\in\Phi^-$, that is, $\delta \in N(ys)$.
At last, we must have $\delta \operatorname{dom} s\gamma$.
Otherwise, if $s \gamma \operatorname{dom} \delta$, then the fact $\beta \operatorname{dom} s \gamma$ gives $\beta \operatorname{dom} \delta$, which implies $\operatorname{B}(\beta,\delta)\ge1$ by Lemma \ref{Bge1}(2), contradicting the fact that $\mathcal{Q}_z$ is 2-spherical.

Suppose that $i$ is odd.
By induction, we have $\gamma \operatorname{dom} s\beta$ and $\operatorname{B}(\gamma,\alpha_s)>0$, and so $\operatorname{B}(\delta,\alpha_s)<0$.
If $\delta \operatorname{dom} s\gamma$, then $\delta \operatorname{dom} \beta$, and then $\operatorname{B}(\delta, \beta) \ge 1$, contradicting the fact that $\beta,\delta\in\mathcal{Q}_z$ again.
Therefore, $s \gamma \operatorname{dom} \delta$ (hence $\gamma \operatorname{dom} s \delta$).
Together with the fact $y\gamma\in\Phi^-$, this in turn gives $ys\delta\in\Phi^-$.
The proof of Step 3 is done.

Combining the above steps, we complete the proof of (iii), and hence complete the proof of Theorem \ref{exist-2-sphere-set}.

\bibliographystyle{amsplain}
\bibliography{boundedness}

@book {And05,
    AUTHOR = {Anderson, James W.},
     TITLE = {Hyperbolic geometry},
    SERIES = {Springer Undergraduate Mathematics Series},
   EDITION = {Second},
 PUBLISHER = {Springer-Verlag London, Ltd.}, 
  address = {London},
      YEAR = {2005},
     PAGES = {xii+276}}

@article {Belolipetsky,
    AUTHOR = {Belolipetsky, Mikhail},
     TITLE = {Cells and representations of right-angled {C}oxeter groups},
   JOURNAL = {Selecta Math.\ (N.S.)},
  FJOURNAL = {Selecta Mathematica. New Series},
    VOLUME = {10},
      YEAR = {2004},
    NUMBER = {3},
     PAGES = {325--339}}

@book{BB05,
 author = {Bj{\"o}rner, Anders and Brenti, Francesco},
 title = {Combinatorics of {Coxeter} groups},
 fseries = {Graduate Texts in Mathematics},
 series = {Grad.\ Texts Math.},
 volume = {231},
 year = {2005},
 publisher = {Springer},
 address = {New York}}

@article {BH93,
    AUTHOR = {Brink, Brigitte and Howlett, Robert B.\},
     TITLE = {A finiteness property and an automatic structure for {C}oxeter
              groups},
   JOURNAL = {Math.\ Ann.},
  FJOURNAL = {Mathematische Annalen},
    VOLUME = {296},
      YEAR = {1993},
    NUMBER = {1},
     PAGES = {179--190}}

@book {Bourbaki-Lie456,
    AUTHOR = {Bourbaki, Nicolas},
    TITLE = {{Lie Groups and Lie Algebras. Chapters 4--6}},
    SERIES = {Elements of Mathematics},
    NOTE = {Translated from the 1968 French original by Andrew Pressley},
    PUBLISHER = {Springer-Verlag},
    address={Berlin},
    YEAR = {2002},
    PAGES = {xii+300}}

@article {EW,
    AUTHOR = {Elias, Ben and Williamson, Geordie},
     TITLE = {The {H}odge theory of {S}oergel bimodules},
   JOURNAL = {Ann.\ of Math.\ (2)},
  FJOURNAL = {Annals of Mathematics. Second Series},
    VOLUME = {180},
      YEAR = {2014},
    NUMBER = {3},
     PAGES = {1089--1136}}

@article{Fu12,
    Author = {Fu, Xiang},
    TITLE = {The dominance hierarchy in root systems of {C}oxeter groups},
    JOURNAL = {J.\ algebra},
    FJOURNAL = {Journal of Algebra},
    VOLUME = {366},
    YEAR = {2012},
    PAGES = {187--204}}

@article {Gao22,
    AUTHOR = {Gao, Jianwei},
     TITLE = {The boundness and lowest two-sided cell of weighted {C}oxeter
              groups of rank 3},
   JOURNAL = {J.\ Lie Theory},
  FJOURNAL = {Journal of Lie Theory},
    VOLUME = {32},
      YEAR = {2022},
    NUMBER = {2},
     PAGES = {499--518}}

@article {Geck11,
    AUTHOR = {Geck, Meinolf},
     TITLE = {On {I}wahori-{H}ecke algebras with unequal parameters and
              {L}usztig's isomorphism theorem},
   JOURNAL = {Pure Appl.\ Math.\ Q.},
  FJOURNAL = {Pure and Applied Mathematics Quarterly},
    VOLUME = {7},
      YEAR = {2011},
    NUMBER = {3},
     PAGES = {587--620}}

@article {GX21,
    AUTHOR = {Gao, Jianwei and Xie, Xun},
     TITLE = {Conjectures {P}1--{P}15 for hyperbolic {C}oxeter groups of
              rank 3},
   JOURNAL = {J.\ Algebra},
  FJOURNAL = {Journal of Algebra},
    VOLUME = {567},
      YEAR = {2021},
     PAGES = {139--195}}

@book {GRS-ramsey,
    AUTHOR = {Graham, Ronald L. and Rothschild, Bruce L. and Spencer, Joel
              H.},
     TITLE = {Ramsey theory},
    SERIES = {Wiley-Interscience Series in Discrete Mathematics and
              Optimization},
   EDITION = {Second},
      NOTE = {A Wiley-Interscience Publication},
 PUBLISHER = {John Wiley \& Sons, Inc.},
   address = {New York},
      YEAR = {1990},
     PAGES = {xii+196}}

@article {GP19-G2,
    AUTHOR = {Guilhot, J\'er\'emie and Parkinson, James},
     TITLE = {A proof of {L}usztig's conjectures for affine type {$G_2$}
              with arbitrary parameters},
   JOURNAL = {Proc.\ Lond.\ Math.\ Soc.\ (3)},
  FJOURNAL = {Proceedings of the London Mathematical Society. Third Series},
    VOLUME = {118},
      YEAR = {2019},
    NUMBER = {5},
     PAGES = {1188--1244}}

@unpublished{HHS,
    author = {Hu, Hongsheng},
    title = {Representations of {C}oxeter groups of {L}usztig's {$\boldsymbol{a}$}-function value $1$},
    note = {to appear, Comm.\ Algebra},
    year = {2026}}

@book {Hum,
    AUTHOR = {Humphreys, James E.},
     TITLE = {Reflection groups and {C}oxeter groups},
    SERIES = {Cambridge Studies in Advanced Mathematics},
    VOLUME = {29},
 PUBLISHER = {Cambridge University Press}, 
   address = {Cambridge},
      YEAR = {1990},
     PAGES = {xii+204}}

@article {Kr,
    AUTHOR = {Krammer, Daan},
     TITLE = {The conjugacy problem for {C}oxeter groups},
   JOURNAL = {Groups Geom.\ Dyn.},
  FJOURNAL = {Groups, Geometry, and Dynamics},
    VOLUME = {3},
      YEAR = {2009},
    NUMBER = {1},
     PAGES = {71--171}}

@incollection {Lus85,
    AUTHOR = {Lusztig, George},
     TITLE = {Cells in affine {W}eyl groups},
 BOOKTITLE = {Algebraic groups and related topics ({K}yoto/{N}agoya, 1983)},
    SERIES = {Adv.\ Stud.\ Pure Math.},
    VOLUME = {6},
     PAGES = {255--287},
 PUBLISHER = {North-Holland}, 
   address = {Amsterdam},
      YEAR = {1985}}

@article {Lus87,
    AUTHOR = {Lusztig, George},
     TITLE = {Cells in affine {W}eyl groups. {II}},
   JOURNAL = {J.\ Algebra},
  FJOURNAL = {Journal of Algebra},
    VOLUME = {109},
      YEAR = {1987},
    NUMBER = {2},
     PAGES = {536--548}}

@book {Lus03,
    AUTHOR = {Lusztig, George},
     TITLE = {Hecke algebras with unequal parameters},
    SERIES = {CRM Monograph Series},
    VOLUME = {18},
 PUBLISHER = {American Mathematical Society}, 
   address = {Providence, RI},
      YEAR = {2003},
     PAGES = {vi+136},
      note = {updated version: arXiv:math/0208154v2, 2014}}

@unpublished{Lus20,
    author = {Lusztig, George},
    title = {Open problems on {I}wahori--{H}ecke algebras},
    note = {preprint, arXiv:2006.08535},
    year = {2020}}

@article {LS,
    AUTHOR = {Li, Yan and Shi, Jianyi},
     TITLE = {The boundedness of a weighted {C}oxeter group with
              non-3-edge-labeling graph},
   JOURNAL = {J.\ Algebra Appl.},
  FJOURNAL = {Journal of Algebra and its Applications},
    VOLUME = {18},
      YEAR = {2019},
    NUMBER = {5},
     PAGES = {1950085, 43}}

@article {SY15,
    AUTHOR = {Shi, Jianyi and Yang, Gao},
     TITLE = {The weighted universal {C}oxeter group and some related
              conjectures of {L}usztig},
   JOURNAL = {J.\ Algebra},
  FJOURNAL = {Journal of Algebra},
    VOLUME = {441},
      YEAR = {2015},
     PAGES = {678--694}}

@article {SY,
    AUTHOR = {Shi, Jianyi and Yang, Gao},
     TITLE = {The boundedness of the weighted {C}oxeter group with complete
              graph},
   JOURNAL = {Proc.\ Amer.\ Math.\ Soc.},
  FJOURNAL = {Proceedings of the American Mathematical Society},
    VOLUME = {144},
      YEAR = {2016},
    NUMBER = {11},
     PAGES = {4573--4581}}

@article {Xi12,
    AUTHOR = {Xi, Nanhua},
     TITLE = {Lusztig's {$A$}-function for {C}oxeter groups with complete
              graphs},
   JOURNAL = {Bull.\ Inst.\ Math.\ Acad.\ Sin.\ (N.S.)},
  FJOURNAL = {Bulletin of the Institute of Mathematics. Academia Sinica. New
              Series},
    VOLUME = {7},
      YEAR = {2012},
    NUMBER = {1},
     PAGES = {71--90}}

@book {Xi94,
    AUTHOR = {Xi, Nanhua},
     TITLE = {Representations of affine {H}ecke algebras},
    SERIES = {Lecture Notes in Mathematics},
    VOLUME = {1587},
 PUBLISHER = {Springer-Verlag},
   address = {Berlin},
      YEAR = {1994},
     PAGES = {viii+137}}

@article {Xie17,
    AUTHOR = {Xie, Xun},
     TITLE = {The lowest two-sided cell of a {C}oxeter group with complete
              graph},
   JOURNAL = {J.\ Algebra},
  FJOURNAL = {Journal of Algebra},
    VOLUME = {489},
      YEAR = {2017},
     PAGES = {38--58}}

@article {Xie21,
    AUTHOR = {Xie, Xun},
     TITLE = {Conjectures {P}1-{P}15 for {C}oxeter groups with complete
              graph},
   JOURNAL = {Adv.\ Math.},
  FJOURNAL = {Advances in Mathematics},
    VOLUME = {379},
      YEAR = {2021},
     PAGES = {Paper No.\ 107565, 62}}

@article {Z,
    AUTHOR = {Zhou, Peipei},
     TITLE = {Lusztig's {$a$}-function for {C}oxeter groups of rank 3},
   JOURNAL = {J.\ Algebra},
  FJOURNAL = {Journal of Algebra},
    VOLUME = {384},
      YEAR = {2013},
     PAGES = {169--193}}

\end{document}